\documentclass[11pt]{article}

\usepackage{amssymb}
\usepackage{graphicx}
\usepackage{color}
\usepackage{amsmath,amsthm,amsfonts,epsfig,setspace}
\numberwithin{equation}{section}

\newcommand{\ds}{\displaystyle}
\def\nm{\noalign{\medskip}}
\newtheorem{thm}{Theorem}[section]

\newtheorem{definition}{Definition} [section]
\newtheorem{lem}{Lemma}[section]

\newtheorem{prop}{Proposition}[section]

\newtheorem{cond}{Condition}

 \def\p{\partial}
\def \Vh0{\stackrel{\circ}{V}_h} \def\to{\rightarrow}

\def\l{\label}  \def\f{\frac} \def\df{\dfrac} 

   \def\eps{\varepsilon}

\def\l|{\left|}
\def\r|{\right|}

\newcommand{\R}{\mathbb{R}}

\newcommand{\lc}
{\mathrel{\raise2pt\hbox{${\mathop<\limits_{\raise1pt\hbox
{\mbox{$\sim$}}}}$}}}

\newcommand{\gc}
{\mathrel{\raise2pt\hbox{${\mathop>\limits_{\raise1pt\hbox{\mbox{$\sim$}}}}$}}}

\newcommand{\ec}
{\mathrel{\raise2pt\hbox{${\mathop=\limits_{\raise1pt\hbox{\mbox{$\sim$}}}}$}}}

\def\be{\begin{equation}} \def\ee{\end{equation}}

\def\bea{\begin{eqnarray}}  \def\eea{\end{eqnarray}}

\def\beas{\begin{eqnarray*}} \def\eeas{\end{eqnarray*}}

\def\bn{\begin{enumerate}} \def\en{\end{enumerate}}

\def\bd{\begin{description}} \def\ed{\end{description}}

\title{Reconstructing fine details of small objects by using plasmonic spectroscopic data. Part II: The strong interaction regime}
\date{}

\author{
Habib Ammari\thanks{\footnotesize Department of Mathematics, 
ETH Z\"urich, 
R\"amistrasse 101, CH-8092 Z\"urich, Switzerland (habib.ammari@math.ethz.ch, sanghyeon.yu@sam.math.ethz.ch). }
\and   Matias Ruiz\thanks{\footnotesize Department of Mathematics and Applications,
Ecole Normale Sup\'erieure, 45 Rue d'Ulm, 75005 Paris, France
(matias.ruiz@ens.fr).}
\and
Sanghyeon Yu \footnotemark[2]
\and
Hai Zhang\thanks{\footnotesize Department of Mathematics, HKUST, Clear Water Bay, Kowloon, Hong Kong (haizhang@ust.hk). The work of Hai Zhang was supported by HK RGC grant ECS 26301016 and startup fund R9355 from HKUST. }
}

\begin{document}
\maketitle

\begin{abstract}
This paper is concerned with the inverse problem of reconstructing a small object from far field measurements by using the field interaction with a plasmonic particle which can be viewed as a passive sensor. It is a follow-up of the work [H. Ammari et al., Reconstructing fine details of small objects by using plasmonic spectroscopic data, SIAM J. Imag. Sci., 11 (2018), 1--23],  where the intermediate interaction regime was considered. In that regime, it was shown that the presence of the target object induces  small shifts to the resonant frequencies of the plasmonic particle. These shifts, which can be determined from the far field data, encodes the contracted generalized polarization tensors of the target object, from which one can perform reconstruction beyond the usual resolution limit. The main argument is based on perturbation theory. However, the same argument is no longer applicable in the strong interaction regime as considered in this paper due to the large shift induced by strong field interaction between the particles. We develop a novel technique based on conformal mapping theory to overcome this difficulty. The key is to design a conformal mapping which transforms the two particle system into a shell-core structure, in which the inner dielectric core corresponds to the target object. We show that a perturbation argument can be used to analyze the shift in the resonant frequencies due to the presence of the inner dielectric core. This shift also encodes information of the contracted polarization tensors of the core, from which one can reconstruct its shape, and hence the target object. Our theoretical findings are supplemented by a variety of numerical results based on an efficient optimal control algorithm.  The results of this paper make the mathematical foundation for plasmonic sensing complete.   
  
\end{abstract}

\medskip

\bigskip

\noindent {\footnotesize Mathematics Subject Classification
(MSC2000): 35R30, 35C20.}


\section{Introduction}
The inverse problem of reconstructing fine details of small objects by using far-field measurements is severally ill-posed. There are two fundamental reasons for this: the diffraction limit and the low signal to noise ratio in the measurements. 

Motivated by plasmonic sensing in molecular biology (see \cite{anker} and the references therein), we developed a new methodology to overcome the ill-posedness of this inverse problem in \cite{part1}. The key idea is to use a plasmonic particle to interact with the target object and to propagate its near field information into far-field in terms of the shifts in the plasmonic resonant frequencies. This plasmonic particle can be viewed as a passive sensor in the simplest form.  For such a plasmonic-particle sensor,  one of the most important characterization is the plasmon resonant frequencies associated with it. These resonant frequencies depend not only on the electromagnetic properties of the particle and its size and shape \cite{matias, matias2, kelly, tocho}, but also on the electromagnetic properties of the environment \cite{matias, kelly, link}. It is the last property which enables the sensing application of plasmonic particles. 
%

In \cite{part1}, the target object is modeled by a dielectric particle whose size is much smaller than that of the sensing plamsonic particle. The intermediate regime where the distance of the two particles is comparable to the size of the plasmonic particle was investigated. 
It was shown that the shifts of the plasmonic resonant frequencies of the plasmonic particle is small and a perturbation argument can be used to derive their asymptotic. Based on these asymptotic formulas, one can obtain their explicit dependence on the generalized polarization tensors of the target particle from which one can perform its reconstruction.
However, when the distance between the particles decreases, 
their interactions increases and the induced shifts  
increase in magnitude as well. 
The perturbation argument will cease to work at certain threshold distance, and the characterization for the shifts of resonant frequencies in terms of information of the target particle becomes more complicated. 


In this paper, we aim to extend the above investigation to the strong interaction regime where the distance of between the two particles is comparable to the size of the small particle. In this regime, the near field interactions are strong and the induced large shifts in plasmonic resonant frequencies cannot be analyzed by a perturbation argument. 
In order to overcome this difficulty, we develop a novel technique based on conforming mapping theory. The key is to design a conformal mapping which transforms the two-particle system into a shell-core structure, in which the inner dielectric core corresponds to the target object. We showed that a perturbation argument can be used to analyze the shift in the resonance frequencies due to the presence of the inner dielectric core. This shift also encodes information on the contracted polarization tensors of the core, from which one can reconstruct its shape, and hence the target object.  
The results of this paper make the mathematical foundation for plasmonic sensing complete.  

The conformal mapping technique has been applied to analyze singular plasmonic systems \cite{pendry1,pendry2}. The nearly touching or touching plasmonic particles system exhibit strong field enhancments and shift of the resonances. The inversion mapping which is conformal was  used to transform two circular disks or spheres into more symmetric systems \cite{bonnetier-triki-circles,pendry3,3D-YA}. After the transformation, the problems become easier to solve.
We also refer to \cite{ciraci} for the fundamental limits of the field enhancements. 
For the general-shaped plasmonic particles, the strong shift of the plasmonic resonances was analyzed in \cite{bonnetier-triki-general}.

We remark that the above idea of plasmonic sensing is closely related to that of 
super-resolution in resonant media, where the basic idea is to propagate the near field information into the far field through certain near field coupling with subwavelength resonators. In a recent series of papers \cite{hai, hai2, hai3}, we have shown mathematically how to realize this idea by using weakly coupled subwavelength resonators and achieve super-resolution and super-focusing. The key is that the near field information of sources can be encoded in the subwavelength resonant modes of the system of resonators through the near field coupling. These excited resonant modes can propagate into the far-field and thus makes the super-resolution from far field measurements possible.

This paper is organized as follows. 
In Section \ref{sec-prelim}, we provide basic results on layer potentials and then explain the concept of plasmonic resonances and the (contracted) generalized polarization tensors.
In Section \ref{sec-forward}, we consider the forward scattering problem of the incident field interaction with a system composed of an dielectric particle and a plasmonic particle. We derive the asymptotic of the scattered field in the case of strong regime. In Section \ref{sec-inverse}, we consider the inverse problem of reconstructing the geometry of the dielectric particle. This is done by constructing the contracted generalized polarization tensors of the target particle through the resonance shifts  
induced to the plasmonic particle.  We provide numerical examples to justify our theoretical results and to illustrate the performances of the proposed optimal control reconstruction scheme.

\section{Preliminaries}\label{sec-prelim}

\subsection{Layer potentials}
We denote by $G(x,y)$ the fundamental solution to the Laplacian in the free space $\R^2$, i.e.,
$$G(x,y) = \f{1}{2\pi}\log|x-y|.$$
Let $D$ be a domain $\mathbb{R}^2$  with $\mathcal{C}^{1,\eta}$ boundary for some $\eta>0$, and let $\nu(x)$ be the outward normal for $ x \in \partial D$. 

We define the single layer potential $\mathcal{S}_{D}$ by 
$$
\mathcal{S}_{D} [\varphi](x) =\int_{\p D }  G(x,y)\varphi(y) d\sigma(y) ,   \quad x \in \mathbb{R}^2, 
$$
and the Neumann-Poincar\'{e} (NP) operator $\mathcal{K}_{D}^*$ by:
$$
\mathcal{K}_{D}^* [\varphi](x) = \int_{\p D }  \frac{ \p G }{\p\nu(x)} (x,y) \varphi(y) d\sigma(y) ,   \quad x \in \p D.
$$
The following jump relations hold:
\begin{align}
{\mathcal{S}_D[\varphi]}\big|_+ &= {\mathcal{S}_D[\varphi]}\big|_-,\label{eqn_jump_single1}
\\
\frac{\p\mathcal{S}_D[\varphi]}{\p\nu}\Big|_{\pm} &= (\pm\frac{1}{2} I +\mathcal{K}_D^*)[\varphi]. \label{eqn_jump_single2}
\end{align}
Here, the subscripts $+$ and $-$ indicate the limits from outside and inside $D$, respectively.

Let $H^{1/2}(\p D)$ be the usual Sobolev space and let $H^{-1/2}(\p D)$ be its dual space with respect to the duality pairing $(\cdot,\cdot)_{-\frac{1}{2},\frac{1}{2}}$. We denote by $H^{-1/2}_0(\p D)$ the collection of all $\varphi\in H^{-1/2}(\p D)$ such that $(\varphi,1)_{-\frac{1}{2},\frac{1}{2}}=0$.

The NP operator is bounded from $H^{-1/2}(\p D)$ to $H^{-1/2}(\p D)$.
Moreover, the operator $\lambda I - \mathcal{K}_D^*: L^2(\partial D)
\rightarrow L^2(\partial D)$ is invertible for any $|\lambda| > 1/2$.
Although the NP operator is not self-adjoint on $L^2(\p D)$, it can be symmetrized on $H_0^{-1/2}(\p D)$ with a proper inner product \cite{kang1,matias}. In fact, 
let $\mathcal{H}^*(\p D)$ be the space $H^{-1/2}_0(\p D)$ equipped with the inner product $(\cdot,\cdot)_{\mathcal{H}^*(\p D)}$ defined by
$$
(\varphi,\psi)_{\mathcal{H}^*(\p D)} =  -(\varphi,\mathcal{S}_D[\psi])_{-\frac{1}{2},\frac{1}{2}},
$$
for $\varphi,\psi\in H^{-1/2}(\p D)$.
Then using the Plemelj's symmetrization principle,
$$
\mathcal{S}_D\mathcal{K}_D^*=\mathcal{K}_D\mathcal{S}_D,
$$
it can be shown that the NP operator $\mathcal{K}_D^*$ is self-adjoint in $\mathcal{H}^*$ with the inner product $(\cdot,\cdot)_{\mathcal{H}^*(\p D)}$. It is also known that $\mathcal{K}_D^*$ is compact when the boundary $\p D$ is $C^{1,\eta}$ \cite{kang1}. 
So it admits the following spectral decomposition in $\mathcal{H}^*$ 
\be\label{spectral_decomposition_Kstar}
\mathcal{K}_D^* = \sum_{j=1}^\infty \lambda_j (\cdot,\varphi_j)_{\mathcal{H}^*} \varphi_j,
\ee 
where $\lambda_j$ are the eigenvalues of $\mathcal{K}_D^*$ and $\varphi_j$ are their associated eigenfunctions. Note that the eigenvalues $|\lambda_j|<1/2$ for all $j \geq 1$ .   

\subsection{Electromagnetic scattering in the quasi-static approximation}

Let us consider a particle $D$ embedded in the free space $\mathbb{R}^2$. Equivalently, the particle $D$ in $\mathbb{R}^3$ has a translational symmetry in the direction of $z$-axis.
Let $\epsilon_D$ (and $\eps_m$) be the permittivity of the particle $D$ (and the background), respectively. So the pemittivity distribution $\eps$ is given by
 $$\eps  = \eps_D \chi(D) +  \eps_m\chi(\mathbb{R}^2 \backslash \overline{D}),$$
 where and $\chi(D)$ is  the characteristic function of $D$.
We are interested in the scattering of the electromagnetic fields $(\mathcal{E},\mathcal{H})$ by the particle $D$.

We assume the particle $D$ is small compared to the wavelength of the incident wave. Then we can adopt the quasi-static approximation and  the electromagnetic scattering can be 
described by a scalar function $u$ which is called the  electric potential.  In the vicinity of the particle $D$, the electric field $E$ is approximated as 
$$
E\approx -\nabla u
$$
and the electric potential $u$ satisfies:
\be
    \begin{cases}
        \nabla \cdot (\eps \nabla u) = 0 &\text{ in }\; \mathbb{R}^2, \\[1.5mm]
         u - u^i = O(|x|^{-1}) &\text{ as }\; |x| \rightarrow \infty,
    \end{cases}
    \label{transmission}
\ee
where $u^i$ is the electric potential of a given incident field and satisfies $\Delta u^i=0$ in $\mathbb{R}^2$.
 
  The electric potential $u$ can be represented as (see, for example, \cite{kang1})
\be
    u = u^i+ \mathcal{S}_{D} [\varphi]  \, ,
    \label{scattered}
\ee
where the density $\varphi$ satisfies the boundary integral equation
\be\label{phi_int_equation}
 (\lambda I - \mathcal{K}_D^*)[\varphi] = \frac{\partial u^i}{\partial \nu}\Big|_{\p D}.
\ee
Here, $\lambda$ is given by
 \be \label{deflambda} 
\lambda= \frac{\eps_D+\eps_m}{2(\eps_D-\eps_m)}.\ee

\subsection{Contracted generalized polarization tensors}\label{subsec-CGPT}

In this subsection, we review the concept of the generalized polarization tensors (GPTs).
It is known that the scattered field $u-u^i$ has the following asymptotic expansion in the far-field  \cite[p. 77] {book2}:
\be
    (u - u^i)(x) =  \sum_{|\alpha|, |\beta| \leq 1 } \f{1}{\alpha ! \beta!} \partial^\alpha u^i(0) M_{\alpha \beta}(\lambda, D) \partial^\beta G(x), 
    \quad |x| \rightarrow + \infty,
    \label{scattered2}
\ee
where $M_{\alpha \beta}(\lambda,D)$ is given by
$$M_{\alpha \beta}(\lambda, D) : =  \int_{\partial D} y^\beta (\lambda I - \mathcal{K}_D^*)^{-1} [\frac{\partial x^\alpha}{\partial \nu}](y)\, d\sigma(y), \qquad \alpha, \beta \in \mathbb{N}^2.$$
Here, the coefficient $M_{\alpha \beta}(\lambda, D)$ is called the {\it generalized polarization tensor} \cite{book2}.

Next we consider the simplified version of the GPTs.
For a positive integer
$m$, let $P_m(x)$ be the complex-valued polynomial
\begin{equation}
P_m(x) = (x_1 + ix_2)^m = r^m \cos m\theta +ir^m \sin m\theta, \label{eq:Pdef}
\end{equation}
where we have used the polar coordinates $x = re^{i\theta}$.

We define the {\it contracted generalized polarization tensors }(CGPTs) to be the
following linear combinations of generalized polarization tensors using the polynomials in
\eqref{eq:Pdef}:
\begin{equation} \label{defCGPT}
\begin{array}{l}
M^{cc}_{m,n}(\lambda,D) = \int_{\partial D} \mbox{Re} \{ P_n\}  (\lambda I - \mathcal{K}_D^*)^{-1} [\frac{\partial \,\mbox{Re} \{ P_m\}}{\partial \nu}]\, d\sigma,  \\
M^{cs}_{m,n}(\lambda,D) = \int_{\partial D} \mbox{Im} \{ P_n\}  (\lambda I - \mathcal{K}_D^*)^{-1} [\frac{\partial \,\mbox{Re} \{ P_m\}}{\partial \nu}]\, d\sigma,\\
M^{sc}_{m,n}(\lambda,D) = \int_{\partial D} \mbox{Re} \{ P_n\}  (\lambda I - \mathcal{K}_D^*)^{-1} [\frac{\partial \,\mbox{Im} \{ P_m\}}{\partial \nu}]\, d\sigma,\\
M^{ss}_{m,n}(\lambda,D) = \int_{\partial D} \mbox{Im} \{ P_n\}  (\lambda I - \mathcal{K}_D^*)^{-1} [\frac{\partial \,\mbox{Im} \{ P_m\}}{\partial \nu}]\, d\sigma. \end{array}
\end{equation}
We remark that CGPTs defined above encodes useful information about the shape of the particle $D$ and can be used for its reconstruction. See \cite{book2, gpt1, GPTs, gpt2} for more details.

For convenience, we introduce the following notation. We denote
\beas
M_{m,n}(\lambda,D) &=& \left( \begin{array}{c c}
M_{m,n}^{cc}(\lambda,D) & M_{m,n}^{cs}(\lambda,D)\\
M_{m,n}^{sc}(\lambda,D) & M_{m,n}^{ss}(\lambda,D)
\end{array}\right).
\eeas
It is worth mentioning that the following symmetry holds (see \cite{book2}):
$$
M_{m,n} = M_{n,m}^T.
$$

When $m=n=1$, the matrix $M(\lambda,D):=M_{1,1}(\lambda,D)$ is called the {\it first order polarization tensor}.
Specifically, we have
$$
M(\lambda,D)_{lm} =  \int_{\partial D} y_j (\lambda I - \mathcal{K}_D^*)^{-1} [\nu_i](y)\, d\sigma(y), \quad l,m=1,2.
$$
We also have from \eqref{scattered2} that
$$
(u-u^i)(x) = \frac{ x^T\cdot M(\lambda,D)\nabla u^i}{|x|^2} + O(|x|^{-2}), \quad
\mbox{as }|x|\rightarrow \infty.
$$
So the leading order term in the far-field expansion of the scattered field $u-u^i$ is determined by the
first order polarization tensor $M(\lambda,D)$. The quantity $M(\lambda,D)(-\nabla u^i)$ is called the dipole moment. In fact, the leading order term is the electric potential generated by a point dipole source with dipole moment $M(\lambda,D)(-\nabla u^i)$.


\subsection{Plasmonic resonances}\label{subsec-plasmonic}

Here we explain the plasmonic resonances. We say that the particle $D$ is plasmonic when its permittivity $\eps_D$ has negative real parts. It is known that the permittivity of noble metals, such as gold and silver, has such a property. More precisely, the permittivity $\epsilon_D$ of the plasmonic (or metallic) particle $D$ is often modeled by the following Drude's model:
\be
\eps_D = \eps_D(\omega) = 1-\frac{\omega_p^2}{\omega(\omega+i\gamma)},
\ee
where $\omega$ is the operating frequency. Here, $\omega_p>0$ means the plasma frequency and $\gamma>0$ means the damping parameter. Usually, the parameter $\gamma$ is a very small number. So $\eps_D(\omega)$ also has a small imaginary part.
Note that, when $\omega<\omega_p$, the permittivity $\eps_D$ has a negative real part.
Contrary to plasmonic particles, ordinary dielectric particles have positive real parts.
Note that, by \eqref{deflambda}, $\lambda$ becomes frequency dependent.


Now we discuss the resonant behavior of the solution $u$ when $\epsilon_D$ is negative (or the particle $D$ is plasmonic).
Recall that  the solution $u$ is represented as
\be
    u = u^i+ \mathcal{S}_{D} [\varphi],  
\ee
where the density $\varphi$ satisfies the boundary integral equation
\be
 (\lambda(\omega) I - \mathcal{K}_D^*)[\varphi] = \frac{\partial u^i}{\partial \nu}\Big|_{\p D}.
\ee
By the spectral decomposition \eqref{spectral_decomposition_Kstar} of $\mathcal{K}_D^*$, we have from \eqref{phi_int_equation} that
 \be\label{varphi_spectral_decomposition}
u = u^i + \sum_{j=1}^\infty \frac{(\frac{\p u^i}{\p\nu},\varphi_j)_{\mathcal{H}^*(\p D)}}{\lambda(\omega)- \lambda_j} \mathcal{S}_D[\varphi_j].
\ee
Recall that $\lambda_j$ are eigenvalues of $\mathcal{K}_D^*$ and they satisfy the condition that $|\lambda_j|<1/2$. 
When $\epsilon_D$ has  negative real parts, we have $|\mbox{Re}\{\lambda(\omega)\}|<1/2$. 
Let $\omega_j$ be such that $\lambda(\omega_j)= \lambda_j$.
Then, if $\omega$ is close to $\omega_j$ and  $(\frac{\p u^i}{\p\nu},\varphi_j)_{\mathcal{H}^*(\p D)}\neq 0$, the function $\mathcal{S}_D[\varphi_j]$ in \eqref{varphi_spectral_decomposition} will be greatly amplified and dominates over other terms. As a result, the magnitude of the scattered field $u-u^i$ will show a pronounced peak at the frequency $\omega_j$ as a function of the frequency $\omega$. This phenomenon is called the plasmonic resonance and $\omega_j$ is called the plasmonic resonant frequency and  $\mathcal{S}_D[\varphi_j]$ is called the resonant mode. 

Let us discuss how we can measure the resonant frequency $\omega_j$ or the eigenvalue $\lambda_j$ from the far field measurements. In fact, the far field for the solution $-\nabla u$ is not equal to the true far-field of the electromagnetic wave since the quasi-static approximation is valid only in the vicinity of the particle $D$. 
But, the polarization tensor $M(\lambda,D)$, which is introduced in the quasi-static approximation, is useful when describing the far-field behavior of the true scattered field.

We first represent $M(\lambda,D)$ in a spectral form.
By \eqref{spectral_decomposition_Kstar}, we have
$$
M(\lambda,D)_{lm} = \sum_{j=1}^\infty \frac{(y_m,\varphi_j)_{-\frac{1}{2},\frac{1}{2}}(\varphi_j,\nu_l)_{\mathcal{H}^*(\p D)}}{\lambda(\omega)-\lambda_j}.
$$
As discussed in Subsection \ref{subsec-CGPT}, the small particle $D$ can be considered as a point dipole source located at $x_0\in \mathbb{R}^2$ and its dipole moment is given by $p_D = M(\lambda,D)(-\nabla u^i)$. We can see from the above spectral representation that the dipole moment $p_D$ becomes resonant when $\omega\approx\omega_j$.

Let $\mathcal{G}^\omega$ be the dyadic Green's function 
$$
\mathcal{G}^\omega(x,y) = (\omega^2 I + \nabla \cdot \nabla ) G^\omega(x,y) 
$$
where $ G^\omega(x,y) = -\frac{i}{4} H_0^{(1)}(\omega |x-y|)$. Then the (true) scattered electric field $E^s$ is well approximated over the whole region as \cite{khelifi, darko}
$$
E^s\approx \mathcal{G}^\omega (x,x_0) p_D.
$$
So, if $\omega\approx \omega_j$, then the amplitude of the scattered wave $E^s$ will be greatly enhanced. So, as a function of the frequency $\omega$, it will have local peaks from which we can recover the resonant frequency $\omega_j$ (or the plasmonic eigenvalue  $\lambda_j$). More specifically, we measure the so called the absorption cross section $\sigma_a$ from the scattered field $E^s$ at the far field region. In fact, this quantity can be approximated as $\sigma_a \propto {\rm{Im}}(p_D)$ for a small plasmonic particle.


\section{The forward problem} \label{sec-forward}

We consider a system composed of a dielectric particle and a plasmonic particle embedded in a homogeneous medium.  
The target dielectric particle and the plasmonic particle  occupy respectively a bounded and simply connected domain $D_1\subset\mathbb{R}^2$ and $D_2\subset\mathbb{R}^2$ of class $\mathcal{C}^{1,\alpha}$ for some $0<\alpha<1$. 
We denote the permittivity of the dielectric particle $D_1$ and the plasmonic particle $D_2$ by $\eps_1$ and $\eps_2$, respectively.
As before, the permittivity of the background medium is denoted by $\varepsilon_m$.
So the permittivity distribution $\eps$ is given by
$$
\eps:=
\eps_1\chi(D_1)  + \eps_2\chi(D_2) +  \eps_m\chi(\R^2\backslash( \overline{D_1\cup D_2})).
$$
As in Subsection \ref{subsec-plasmonic}, the permittivity $\eps_2$ of the plasmonic particle $D_2$ depends on the operating frequency and is modeled as
$$
\eps_2=\eps_2(\omega) = 1-\frac{\omega_p^2}{\omega(\omega+i\gamma)}.
$$

The total electric potential $u$ satisfies the following equation:
\be \label{eq-Helmholtz eq biosensing}
\begin{cases}
	\ds \nabla \cdot (\eps \nabla u)  = 0 &\quad \mbox{in } \R^2\backslash (\p D_1\cup \p D_2), \\
	\nm
	 u|_{+} = u|_{-}    &\quad \mbox{on } \partial D_1 \cup \partial D_2, \\
	\nm
	  \ds \eps_{m} \df{\p u}{\p \nu} \Big|_{+} =\eps_{1} \df{\p u}{\p \nu} \Big|_{-}  &\quad \mbox{on } \partial D_1,\\
	\nm
	  \ds \eps_{m} \df{\p u}{\p \nu} \Big|_{+} = \eps_{2} \df{\p u}{\p \nu} \Big|_{-}  &\quad \mbox{on } \partial D_2,\\
	\nm
	  (u-u^i)(x) = O(|x|^{-1}),   &\quad\mbox{as }|x|\rightarrow \infty,
\end{cases}
\ee
where $u^i(x)$ is the electric potential for a given incident field as before.

\subsection{Boundary integral formulation}
We derive a layer potential representation of the total field $u$ to \eqref{eq-Helmholtz eq biosensing} in this section. We first denote by 
$u_{D_1}$ the total field resulting from the incident field $u^i$ and the ordinary particle $D_1$ 
(in the absence of the plasmonic particle $D_2$).
Let us denote 
$$
\lambda_{D_j} = \df{\eps_j+\eps_m}{2(\eps_j-\eps_m)}, \quad j=1, 2.
$$
Then $u_{D_1}$ has the following representation \cite{book2}: 
\beas
u_{D_1}(x) = u^i(x) + \mathcal{S}_{D_1}\left(\lambda_{D_1}Id - \mathcal{K}_{D_1}^*\right)^{-1}[\df{\p u^i}{\p \nu_1}](x), \quad \mbox{ for }  x \in \R^2 \backslash \overline{D_1}.
\eeas

We next introduce the Green function $G_{D_1}(\cdot,y)$ for the medium with 
permittivity distribution $ \eps_{D_1} \chi(D_1) +  \eps_m\chi(\mathbb{R}^2 \backslash \overline{D_1})$. More precisely, $G_{D_1}(\cdot,y)$ satisfies the following equation
\[
\nabla_x \cdot \left((\eps_{D_1} \chi(D_1) +  \eps_m\chi(\mathbb{R}^2 \backslash \overline{D_1})) \nabla_x G_{D_1}(x, y)\right)  = \delta(x-y).
\]
Using $G_{D_1}$, we define the layer potential $\mathcal{S}_{D_2,D_1}$ by 
\beas
\mathcal{S}_{D_2,D_1}[\varphi](x) = \int_{\p D_2}G_{D_1}(x,y)\varphi(y)d\sigma(y).
\eeas
We also define 
$$
\mathcal{A}= \mathcal{K}_{D_2}^* - \f{\p }{\p \nu_2}\mathcal{S}_{D_1}\left(\lambda_{D_1}Id - \mathcal{K}_{D_1}^*\right)^{-1}\df{\p \mathcal{S}_{D_2}[\cdot]}{\p \nu_1}.
$$

It was proved in \cite{part1} that the solution $u$ can be represented using $\mathcal{S}_{D_2,D_1}$ and $\mathcal{A}$ as shown in the following lemma.

\begin{lem}\cite{part1}
The total electric potential $u$ can be represented as follows:
\be \label{solution helm nanoparticle biosenging}
u = u_{D_1} + \mathcal{S}_{D_2,D_1} [\psi], \quad x \in \R^2 \backslash \overline{D_2},
\ee
where the density $\psi$ satisfies 

\be\label{psi_integral_equation2}
\left(\lambda_{D_2}Id-\mathcal{A}\right) [\psi] = \df{\p u_{D_1}}{\p \nu_2}.
\ee

\end{lem}

\subsection{Strong interaction regime and conformal transformation}\label{sec-strong-conformal}

We assume the following condition on the sizes of the particles $D_1$ and $D_2$.

\begin{cond} \label{condition0 biosensing} 
 The plasmonic particle $D_2$ has size of order one;  the dielectric particle $D_1$ has size of order $\delta \ll 1$.
 \end{cond}

\begin{definition}[\textbf{Strong interaction regime}] \label{def-strong sensing}
We say that the small dielectric particle $D_1$ is in the strong regime with respect to the plasmonic particle $D_2$ if there exist positive constants $C_1$ and $C_2$ such that $C_1 < C_2$ and
$$
C_1 \delta \leq 
{\rm{dist}} (D_1,D_2) \leq C_2 \delta.
$$
\end{definition}

Definition \ref{def-strong sensing} says that the dielectric particle $D_1$ is closely located  to the plasmonic particle $D_2$ with a separation distance of order $\delta$.

In our recent paper \cite{part1}, the intermediate interaction regime is considered. The key observation is that, if we assume the distance between $D_1$ and $D_2$ is assumed to be of order one, then the effect of the small unknown particle $D_1$ can be considered as a small perturbation. To see this, we rewrite the equation (\ref{psi_integral_equation2}) in the form 
\be\label{psi_integral_equation22}
\left(\mathcal{A}_{D_2,0} + \mathcal{A}_{D_2,1}\right)[\psi] = \df{\p u_{D_1}}{\p \nu_2},
\ee
where
\bea
\mathcal{A}_{D_2,0} &=& \lambda_{D_2}Id - \mathcal{K}_{D_2}^*,\nonumber\\
\mathcal{A}_{D_2,1} &=&\: \f{\p }{\p \nu_2}\mathcal{S}_{D_1}\left(\lambda_{D_1}Id - \mathcal{K}_{D_1}^*\right)^{-1}\df{\p \mathcal{S}_{D_2}[\cdot]}{\p \nu_1}. \label{def-perturbation operator biosensing}
\eea
It can be shown that the operator $\mathcal{A}_{D_2,1}$ 
is a small perturbation to the operator $\mathcal{A}_{D_2,0}$ \cite{part1}, and so, the authors were able to apply the perturbation method for analyzing the plasmonic resonance.
However, in the strong interaction regime, the operator $\mathcal{A}_{D_2,1}$ is no longer small compared to the latter. 
As a consequence, the perturbation theory is not applicable and it becomes challenging to analyze the interaction between the particles .

We now introduce a method to tackle this issue by using conformal mapping technique.
Let $B_1$ be a circular disk containing the dielectric particle $D_1$ with radius $r_1$ of order $\delta$.
We assume the plasmonic particle $D_2$ is a circular disk with radius $r_2$.
For convenience, we denote $D_2$ by $B_2$.
We emphasize that the shape of $D_1$ is unknown.
We let $d$ to be the distance between the two disks $B_1$ and $B_2$, {\it i.e.},
$$
d=\mbox{dist}(B_1,B_2).
$$
By the assumption, $d$ is of order $\delta$.

Let $R_j$ be the reflection with respect to $\p B_j$ and let $\mathbf{p}_1$ and $\mathbf{p}_2$ be the unique fixed points of
the combined reflections $R_1\circ R_2$ and $R_2\circ R_1$, respectively.
Let $\mathbf{n}$ be the unit vector in the direction of $\mathbf{p}_2-\mathbf{p}_1$.
We set $(x,y)\in\mathbb{R}^2$ to be the Cartesian coordinates such that
$\mathbf{p}=(\mathbf{p}_1+\mathbf{p}_2)/2$ is the origin and the $x$-axis is parallel to $\mathbf{n}$.
Then one can see that $\mathbf{p}_1$ and $\mathbf{p}_2$ can be written as
\be\label{pjdef}
\mathbf{p}_1=(-a,0)\quad\mbox{and}\quad\mathbf{p}_2=(a,0),
\ee
where the constant $a$ is given by
\be\label{a_def}
a =\frac{\sqrt d\sqrt{ (2 r_1 + d) (2 r_2 + d) (2 r_1 + 2 r_2 +
  d)}}{2 (r_1 + r_2 + d)}.
\ee
Then the center $\mathbf{c}_i$ of $B_i$ ($i=1,2$) is given by
\be\label{c1c2}
\mathbf{c}_i = \Big((-1)^i \sqrt{r_i^2+a^2},0\Big).
\ee 

Define the conformal transformation $\Phi$ by
$$
\zeta=\Phi(z) = \frac{z+a}{z-a}, \quad z=x+i y.
$$
In other words,
$$
z=\Phi^{-1}(\zeta) = a \frac{\zeta+1}{\zeta-1}.
$$
We also define
$$
s_j = (-1)^j\sinh^{-1}(a/r_j), \quad j=1,2,
$$
and the two disks $\widetilde B_1$ and $\widetilde B_2$ by
$$
\widetilde{B}_1 = \{ |\zeta| < \tilde{r}_j \}, \quad \tilde{r}_j=\exp( s_j), \ j=1,2.
$$
It can be shown that, in the $\zeta$-plane, the disks $B_1$ and $B_2$ are transformed to
$$
\Phi(B_1) = \widetilde{B}_1 = \{ |\zeta| < \tilde{r}_1 \},
$$
and
$$
\Phi(B_2) =\mathbb{R}^2\setminus \overline{\widetilde{B}_2} = \{ |\zeta| > \tilde{r}_2 \}.
$$
One can check that $\tilde{r}_1<1$ and $\tilde{r}_2>1$.
The exterior region $\mathbb{R}^2\setminus\overline{B_1 \cup B_2}$ becomes a shell region between $\p\widetilde{B}_1$ and $\p\widetilde{B}_2$  in the $\zeta$-plane:
$$
\Phi(\mathbb{R}^2\setminus\overline{B_1 \cup B_2}) =\widetilde{B}_2 \setminus\overline{\widetilde{B}_1} =  \{ \tilde{r}_1 < |\zeta| < \tilde{r}_2 \}.
$$

\begin{figure*}
\begin{center}
\epsfig{figure=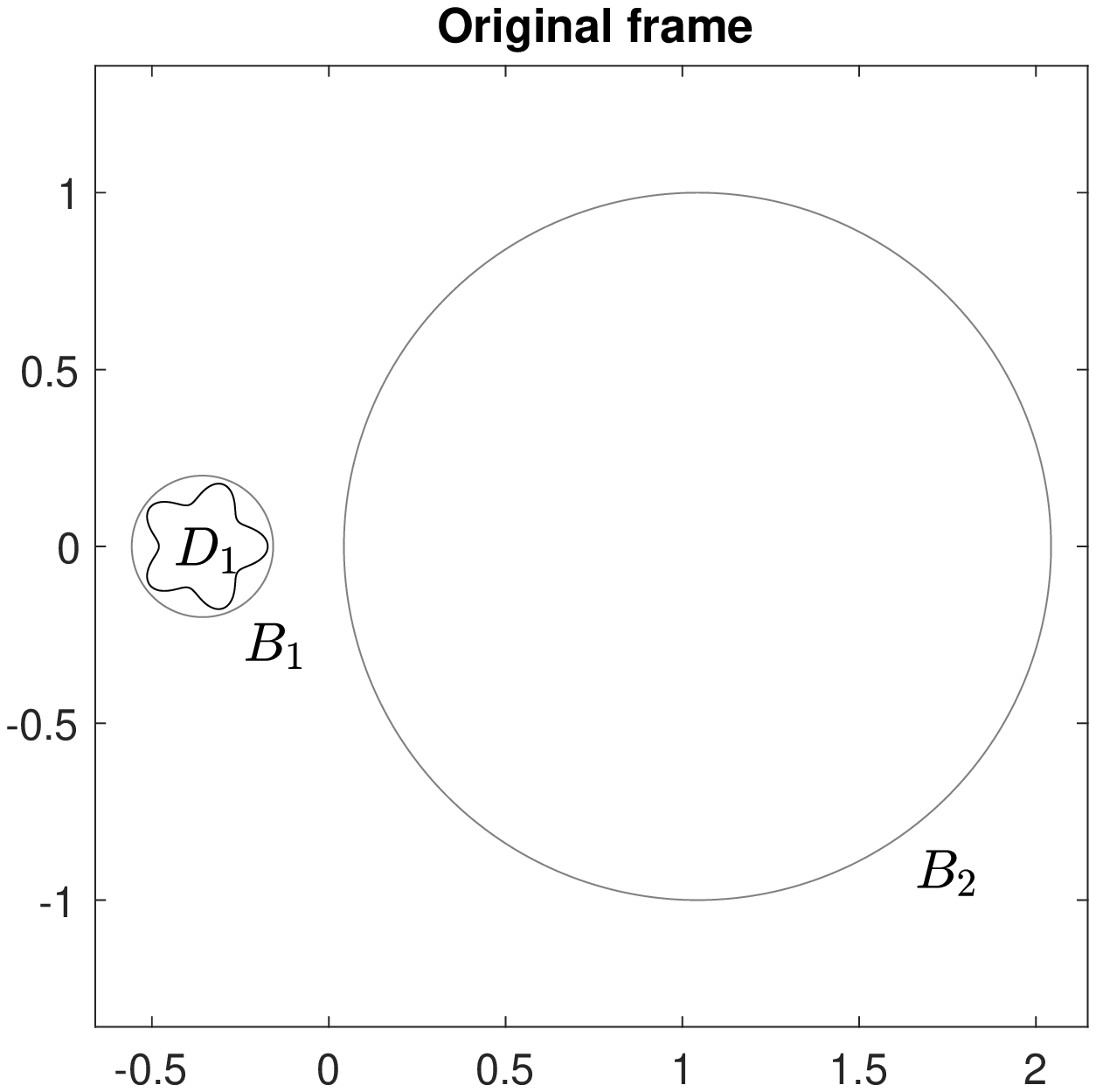,width=8cm}
\hskip-0.5cm
\epsfig{figure=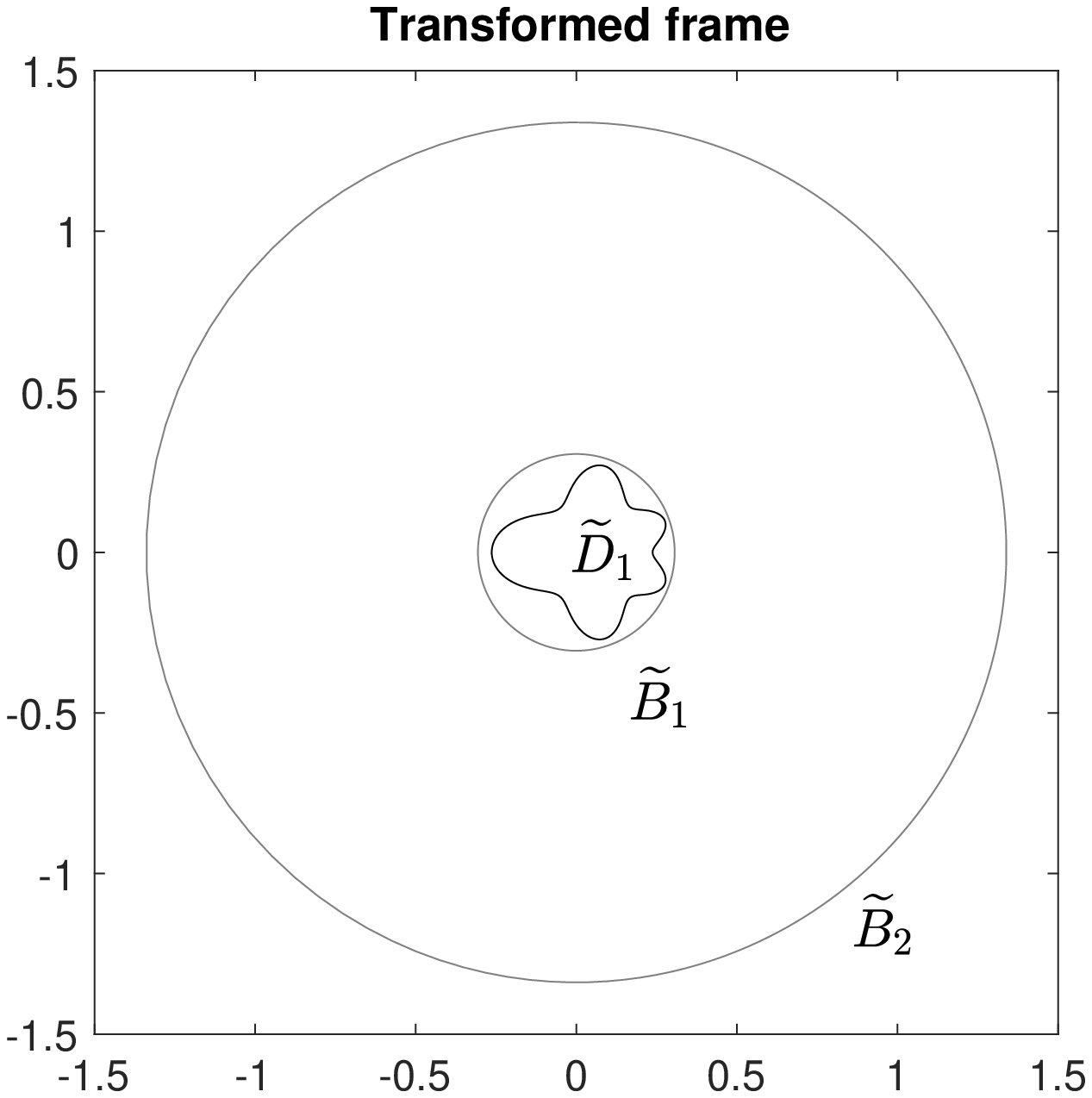,width=8cm}
\end{center}
\caption{ (left) original configuration and (right) transformed one by the conformal map $\Phi$}
\label{fig:config}
\end{figure*}

\begin{figure*}
\begin{center}
\epsfig{figure=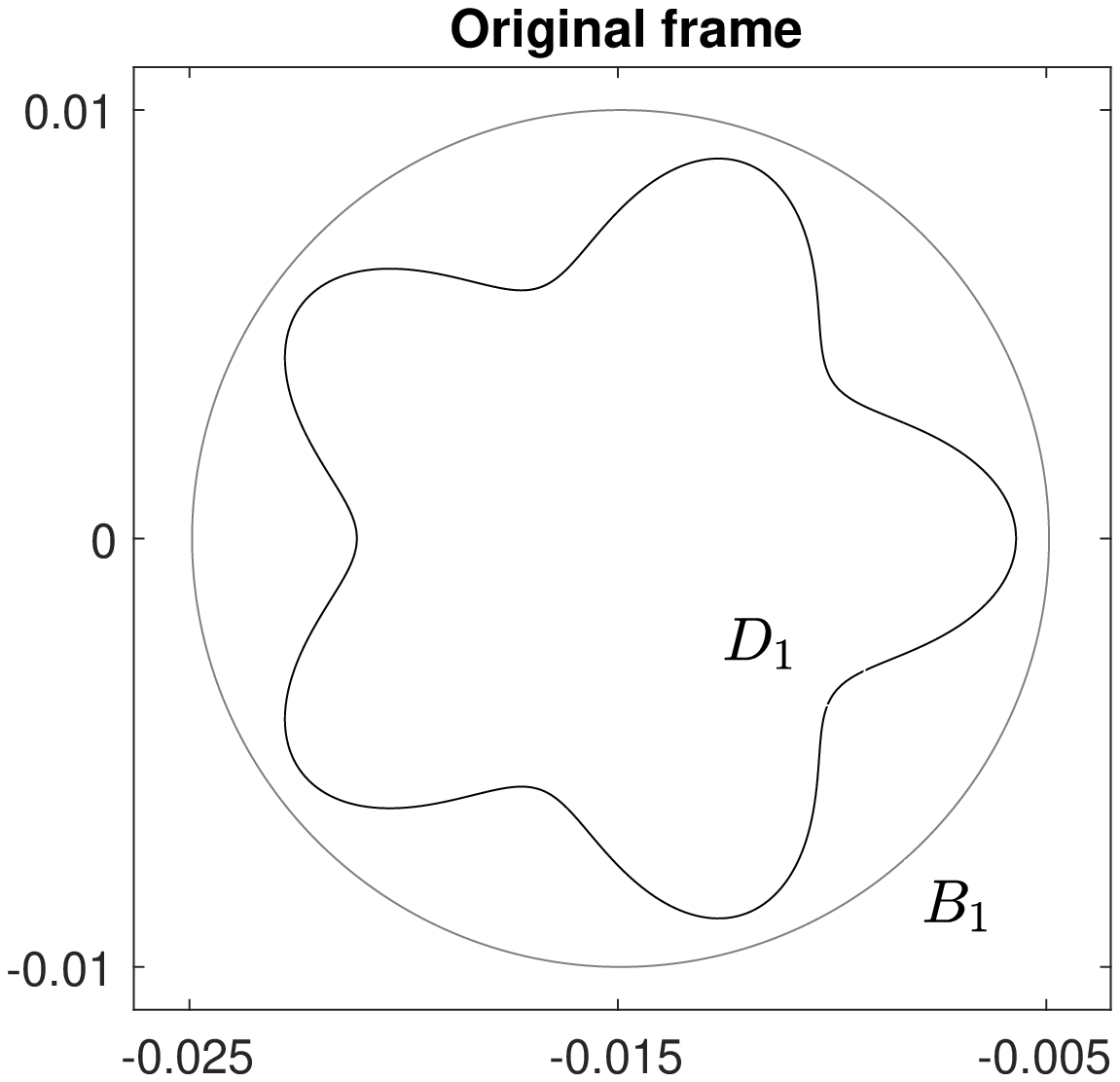,width=5.75cm}
\hskip-0.6cm
\epsfig{figure=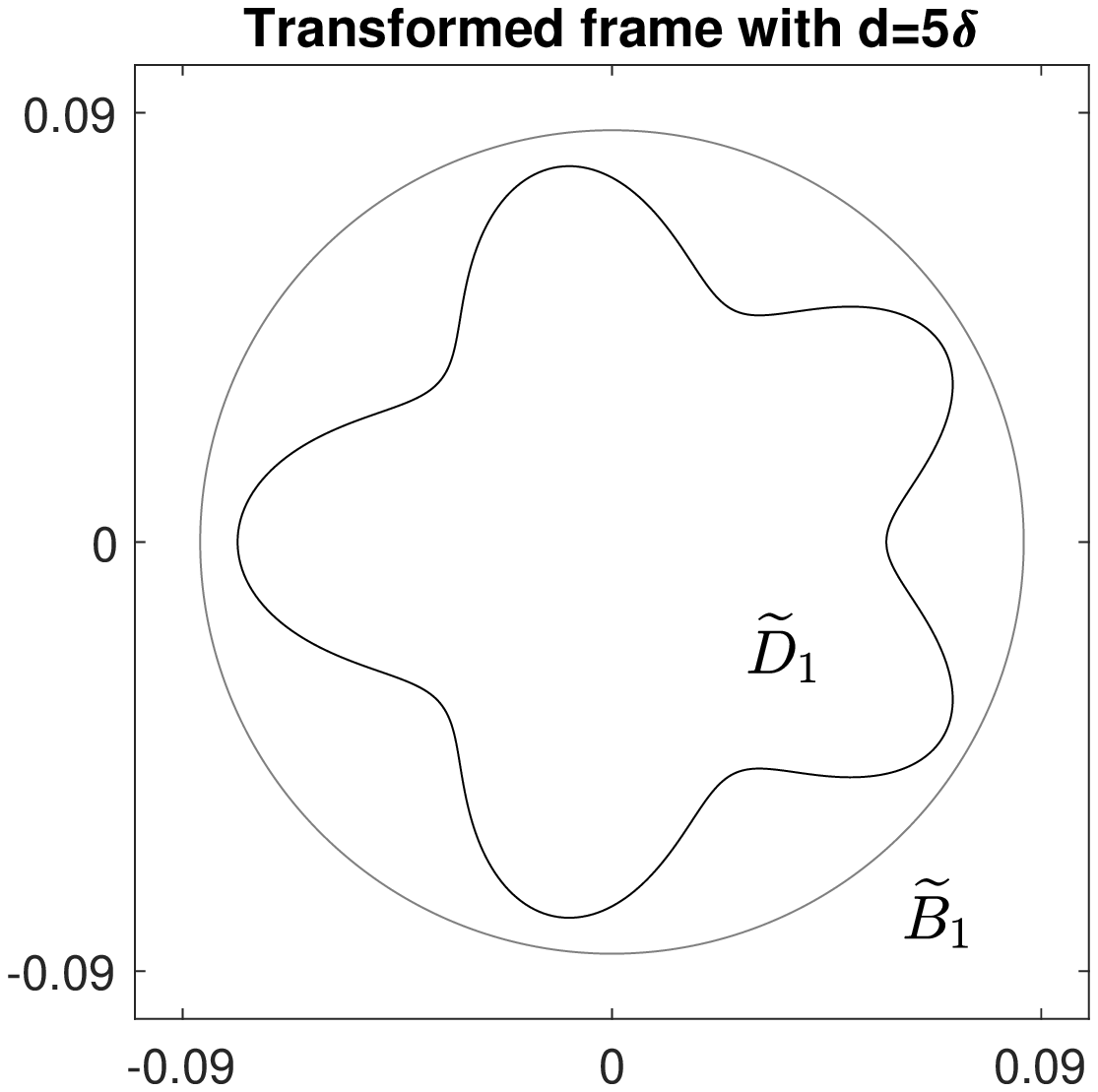,width=5.5cm}
\hskip-0.6cm
\epsfig{figure=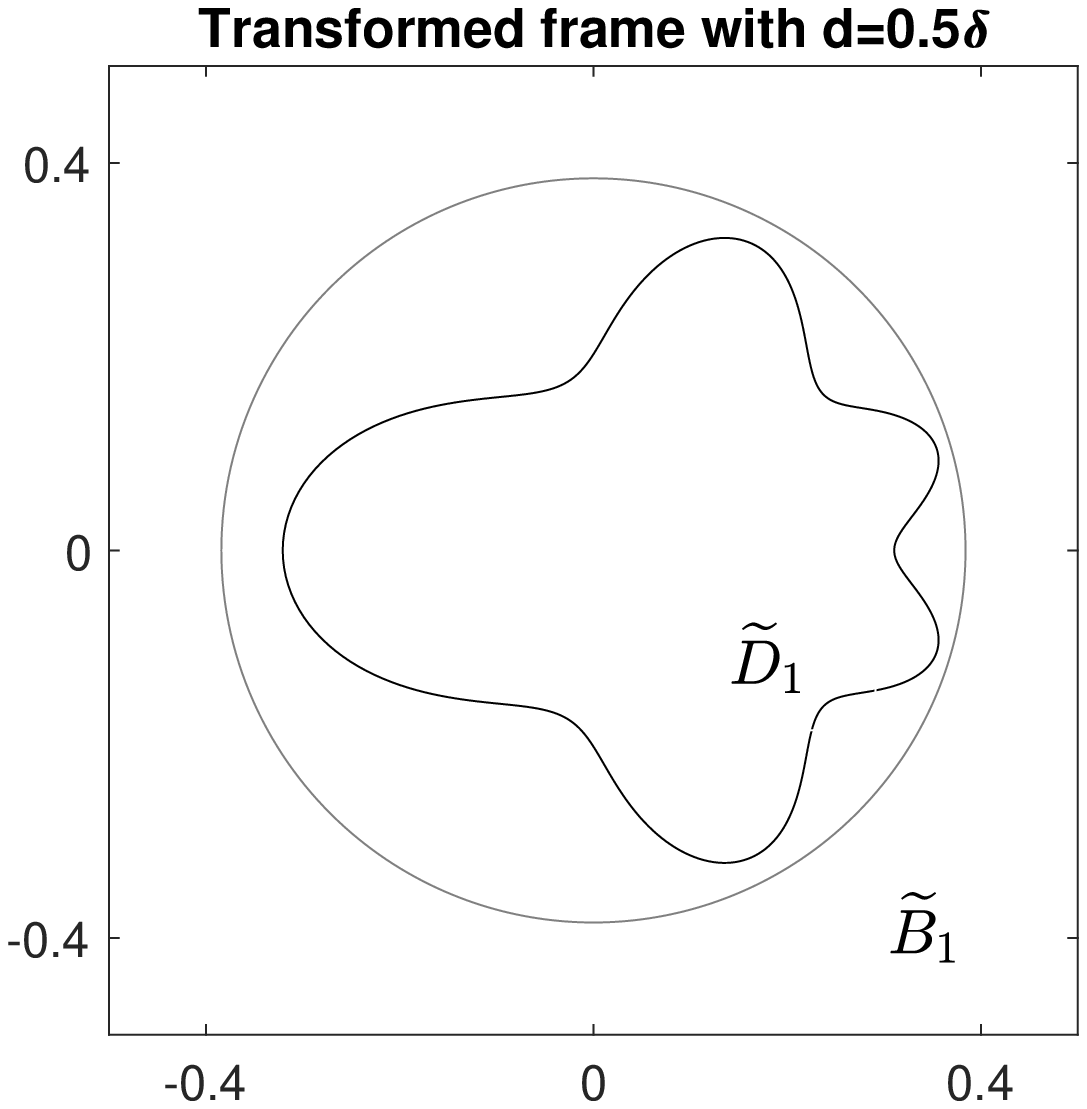,width=5.55cm}
\end{center}
\caption{ (left) original configuration (center) the transformed one with $d=5 \delta$ (right) the same but with $d=0.5 \delta$. We set $r_1 = \delta$, $r_2=1$ and $\delta = 0.01$. }
\label{fig:deform}
\end{figure*}

{
To illustrate the geometry, in Figure \ref{fig:config},we show an example for the configuration of a system of a small dielectric particle $D_1$ and a 
plasmonic particle $B_2$. We also show its transformed geometry by the conformal map $\Phi$. We set $\delta = 0.2$, $r_1 = \delta$, $r_2 = 1$ and $d=\delta$.  

It is worth mentioning that the shape of the transformed domain $\widetilde{D}_1$ strongly depends on the ratio between $d$ and $\delta$ but is independent of $\delta$ itself. 
Suppose that $d=c\delta$ for some $c>0$. If $c$ is of order one, then the shape of $\widetilde{D}_1$ is almost the same as that of $D_1$.
On the contrary, if $c$ is too small, then the shape of $\widetilde{D}_1$ is highly distorted. See Figure \ref{fig:deform}.

}

\subsection{Boundary integral formualtion in the transformed domain}

Let us define $\tilde{u}(\zeta) = u(\Phi^{-1}(\zeta))$ and $\tilde{u}^i(\zeta) = u^i(\Phi^{-1}(\zeta))$. Then, since the mapping $\Phi$ is conformal, $\tilde{u}$ and $\tilde{u}^i$ are harmonic in the $\zeta$-plane. Moreover, the transmission conditions for $\tilde{u}$ are preserved. In fact, the transformed potential $\tilde{u}$ satisfies the following equations:
\be \label{eq-Helmholtz eq biosensing transformed}
\begin{cases}
	\ds \nabla \cdot (\tilde\eps \nabla \tilde{u})  = 0 &\quad \mbox{in } \R^2\backslash (\p \widetilde{D}_1\cup \p \widetilde{D}_2), \\
	\nm
	 \tilde{u}|_{+} = \tilde{u}|_{-}    &\quad \mbox{on } \partial \widetilde{D}_1 \cup \partial \widetilde{D}_2, \\
	\nm
	  \ds \eps_{m} \df{\p \tilde{u}}{\p \nu} \Big|_{+} =\eps_{1} \df{\p \tilde{u}}{\p \nu} \Big|_{-}  &\quad \mbox{on } \partial \widetilde{D}_1,\\
	\nm
	  \ds \eps_{2} \df{\p \tilde{u}}{\p \nu} \Big|_{+} = \eps_{m} \df{\p \tilde{u}}{\p \nu} \Big|_{-}  &\quad \mbox{on } \partial \widetilde{D}_2,\\
	\nm
	  (\tilde{u}-\tilde{u}^i)(\zeta) = O(|\zeta - (1,0)|)   &\quad\mbox{as }\zeta\rightarrow (1,0),
\end{cases}
\ee
where the transformed permittivity distribution $\tilde\eps$ is given by
$$
\tilde\eps=
\eps_1\chi(\widetilde{D}_1)  + \eps_2\chi(\mathbb{R}^2\setminus\widetilde{D}_2) +  \eps_m\chi(\widetilde{D}_2 \setminus\overline{\widetilde{D}_1}).
$$
Note that the transformed problem looks similar to the original one, even though the geometry of the particles is of a completely different nature.
As $\delta$ goes to zero, the radii $\tilde{r}_1$ and $\tilde{r}_2$ have the following asymptotic properties:
$$
\tilde{r}_1 = \tilde{r}_1^0 + O(\delta)
,\quad
\tilde{r}_2 = 1+O(\delta)
$$
for some $0<r_1^0 < 1$ independent of $\delta$.
Hence, in contrast to the original problem, the transformed boundaries $\p \widetilde{B}_1$ and $\p \widetilde{B}_2$ ($=\p\widetilde{D}_2$) are not close to touching. Moreover, they share the same center (see Figure \ref{fig:config}).
This will enable us to analyze more deeply the spectral nature of the problem.

Now we represent the solution to the transformed problem using the layer potentials.
By applying a similar procedure as the one used for \eqref{def-perturbation operator biosensing}, we can obtain the following representation:
\be \label{solution helm nanoparticle biosenging transformed}
\tilde u = (\mbox{const.}) + u_{\widetilde{D}_1} + \mathcal{S}_{\widetilde{D}_2,\widetilde{D}_1} [\widetilde\psi], \quad x \in \R^2.
\ee
Here, the constant term is needed to satisfy the last condition in \eqref{eq-Helmholtz eq biosensing transformed}.
The density function $\widetilde\psi$ satisfies the following boundary integral equation:
\be\label{psi_integral_equation}
\big(\lambda_{D_2}I  - \widetilde{\mathcal{A}}\,\big)\big[\widetilde \psi\big] = \df{\p \tilde u_{\widetilde D_1}}{\p \nu_2},
\ee
where
\bea
\widetilde{\mathcal{A}} &=&  \mathcal{K}_{\widetilde{D}_2}^* - \f{\p }{\p \nu_2}\mathcal{S}_{\widetilde{D}_1}\left(\lambda_{D_1}I - \mathcal{K}_{\widetilde{D}_1}^*\right)^{-1}\df{\p \mathcal{S}_{\widetilde{D}_2}[\cdot]}{\p \nu_1}, \label{def-perturbation operator biosensing2}
\\
\tilde u_{\widetilde{D}_1} &=& \tilde u^i + \mathcal{S}_{\widetilde{D}_1}\left(\lambda_{D_1}I - \mathcal{K}_{\widetilde{D}_1}^*\right)^{-1}\Big[\df{\p \tilde{u}^i}{\p \nu_1}\Big].
\eea

\begin{lem} The following relation between $\mathcal{A}$ and $\widetilde{\mathcal{A}}$ holds
\be\label{A_identical_Atilde}
( \phi,\mathcal{A}[\psi])_{\mathcal{H}^*(\p D_2)} = ( \widetilde\phi,\widetilde{\mathcal{A}}[\widetilde\psi])_{\mathcal{H}^*(\p \widetilde{D}_2)},
\ee
where $\phi,\psi\in \mathcal{H}^*(\p D_2)$ and $\widetilde{\phi}=\phi\circ\Phi^{-1}, \widetilde{\psi}=\psi\circ\Phi^{-1}$.
\end{lem} 
\proof
By the conformality of the map $\Phi$, the single layer potentials $\mathcal{S}_{D_2}[\phi]$ and $\mathcal{S}_{\widetilde{D}_2}[\widetilde\phi]\circ \Phi$ are identical up to an additive constant, whence (\ref{A_identical_Atilde}) follows. 

\subsection{Computation of the operator $\widetilde{\mathcal{A}}$ and its spectral properties}

Here we compute the operator $\widetilde{\mathcal{A}}$. Note that  $\widetilde{\mathcal{A}}$ is an operator which maps $\mathcal{H}^*(\p \widetilde{D}_2)$ onto $\mathcal{H}^*(\p \widetilde{D}_2)$. Since $\p \widetilde{D}_2$ is a circle, we use the Fourier basis for $\mathcal{H}^*(\p \widetilde{D}_2)$. Let $(r,\theta)$ be the polar coordinates in the $\zeta$-plane, {\it i.e.}, $\zeta=r e^{i\theta}$. We define
$$
\varphi_n^c(\theta) = \cos{n\theta}, \quad \varphi_n^s(\theta)=\sin{n\theta}.
$$  

The following proposition holds. 
\begin{prop}\label{prop:comp_Atilde}
We have
\begin{align} \label{form1}
\widetilde{\mathcal{A}} [\varphi_n^c](\zeta) &= \sum_{m=1}^\infty-\frac{\tilde{r}_2^{-(n+m)}}{4\pi n}  (M_{nm}^{cc}(\lambda_{D_1}, \widetilde{D}_1) \cos m\theta + M_{nm}^{cs}(\lambda_{D_1}, \widetilde{D}_1) \sin m\theta),
\end{align}
and
\begin{align} \label{form2}
\widetilde{\mathcal{A}} [\varphi_n^s](\zeta) &= \sum_{m=1}^\infty-\frac{\tilde{r}_2^{-(n+m)}}{4\pi n}  (M_{nm}^{sc}(\lambda_{D_1}, \widetilde{D}_1) \cos m\theta + M_{nm}^{ss}(\lambda_{D_1}, \widetilde{D}_1) \sin m\theta)
\end{align}
for $n\neq 0$.
\end{prop}

\proof
Since $\p \widetilde{D}_2$ is a circle, $\mathcal{K}^*_{\widetilde{D}_2}=0$ on $\mathcal{H}^*(\p\widetilde{D}_2)$.
Therefore, we only need to consider the second term in $\mathcal{A}$.
It is easy to see that
\begin{align}
\mathcal{S}_{\widetilde{D}_2}\big[\varphi_{n}^{c}\big] (r,\theta) 
&=
\ds -\frac{\tilde{r}_2^{-n+1}}{2n}  r^n
{\cos{ n \theta}} , 
\label{aaabbbccc}
\\
\mathcal{S}_{\widetilde{D}_2}\big[\varphi_{n}^{s}\big] (r,\theta) 
&=
\ds -\frac{\tilde{r}_2^{-n+1}}{2n}  r^n
{\sin{ n \theta}} , 
\end{align}
for $0\leq r\leq \tilde{r}_2$. Thus, we have
\begin{align}
\widetilde{\mathcal{A}} [\varphi_n^c](\zeta) &= -\frac{\tilde{r}_2^{-n+1}}{2n}\frac{\p}{\p\nu_2}\int_{\p \widetilde{D}_1 }G(\zeta,\zeta')
\left(\lambda_{D_1}I - \mathcal{K}_{\widetilde{D}_1}^*\right)^{-1}
 \left[\frac{\p}{\p\nu_1}\mbox{Re}\{ P_n \}\right](\zeta') \,d\sigma(\zeta').
\end{align}
It is known that \cite{book3}
$$
G({x},y) = \sum_{m=1}^{\infty} \f{(-1)}{2\pi m}  \f{\cos(m \theta_x)}{ r_x^m} r_{{y}}^m\cos(m \theta_{{y}}) + \f{(-1)}{2\pi m} \f{\sin(m \theta_x)}{ r_x^m} r_{{y}}^m\sin(m \theta_{{y}}), \quad |x|<|y|,
$$
where $(r_x,\theta_x)$ and $(r_y,\theta_y)$ are the polar coordinates of $x$ and $y$, respectively. Then, by letting $x=\zeta$ and $y=\zeta' \in \p \widetilde{D}_2$, we get
\begin{align*}
\widetilde{\mathcal{A}} [\varphi_n^c](\zeta) &= \sum_{m=1}^\infty-\frac{\tilde{r}_2^{-(n+m)}}{4\pi n} \cos m\theta \int_{\p \widetilde{D}_1 }\mbox{Re}\{P_m\}
\left(\lambda_{D_1}I - \mathcal{K}_{\widetilde{D}_1}^*\right)^{-1}
 \left[\frac{\p}{\p\nu_1}\mbox{Re}\{ P_n \}\right](\zeta') \,d\sigma(\zeta')
 \\
 &\quad -\frac{\tilde{r}_2^{-(n+m)}}{4\pi n} \sin m\theta \int_{\p \widetilde{D}_1 }\mbox{Im}\{P_m\}
\left(\lambda_{D_1}I - \mathcal{K}_{\widetilde{D}_1}^*\right)^{-1}
 \left[\frac{\p}{\p\nu_1}\mbox{Re}\{ P_n \}\right](\zeta') \,d\sigma(\zeta').
\end{align*}
Finally, from the definition of the CGPTs (see  (\ref{defCGPT})), (\ref{form1}) follows. Similarly, one can derive (\ref{form2}).
\qed

Let us define
$$
M_{nm}=M_{nm}(\lambda_{D_1} , \widetilde{D}_1) = \begin{pmatrix}
 M_{nm}^{cc}(\lambda_{D_1} , \widetilde{D}_1) & M_{nm}^{cs}(\lambda_{D_1} , \widetilde{D}_1)
 \\
 M_{nm}^{sc}(\lambda_{D_1} , \widetilde{D}_1) & M_{nm}^{ss}(\lambda_{D_1} , \widetilde{D}_1)
 \end{pmatrix},
$$
and
\be\label{tilde_M_nm}
\widetilde{M}_{nm}= -\frac{\tilde{r}_2^{-(n+m)}}{4\pi n} M_{nm}(\lambda_{D_1} , \widetilde{D}_1).
\ee
In view of Proposition \ref{prop:comp_Atilde}, we see that the operator $\widetilde{\mathcal{A}}$ can be represented in a block matrix form as follows:
\begin{equation} 
\widetilde{\mathcal{A}} = \begin{bmatrix}
 \widetilde{M}_{11} & \widetilde{M}_{12} & \widetilde{M}_{13} & \cdots
 \\
 \widetilde{M}_{21} & \widetilde{M}_{22} & \cdots & \cdots
 \\
 \widetilde{M}_{31} & \cdots & \cdots &
 \\
 \cdots & & &
 \end{bmatrix}.
 \label{Atilde_matrix}
\end{equation}

Recall that $\widetilde{D}_1$ is contained in the disk $\widetilde{B}_1$ with radius $\tilde{r}_1$. One can derive that 
$$|M_{nm}(\lambda_{D_1},\widetilde{D}_1)|\leq C \tilde{r}_1^{n+m} $$ for some positive constant $C$  \cite{book2}. Therefore, 
\be\label{Mnm_tilde_decay}
|\widetilde{M}_{nm}(\lambda_{D_1},\widetilde{D}_1)| \leq C \left(\frac{\tilde{r}_1}{\tilde{r}_2}\right)^{n+m}.
\ee
This decay property of $ \widetilde{M}_{nm}$ is crucial for our conformal mapping technique. An important consequence is that
the operator $\widetilde{\mathcal{A}}$ can be efficiently approximated by finite dimensional  matrices obtained through a standard truncation procedure. Here we remark that $\widetilde{\mathcal{A}} = O(({\tilde{r}_1}/{\tilde{r}_2})^2)$.

{
If the particle $D_1$ is in the strong regime, then we may write $d=c\delta$ for some $c>0$. If $c$ is of order one, the ratio $\frac{\tilde{r}_1}{\tilde{r}_2}$ is relatively small (but regardless of how small $\delta$ is). 
 In section \ref{sec-inverse} we apply the eigenvalue perturbation method to analyze the spectral nature more explicitly when we consider the related inverse problem.
}

\subsection{Spectral decomposition $\mathcal{A}$ of and the scattered field }
It is clear that $\widetilde{\mathcal{A}}$ (or $\mathcal{A}$) is compact. Moreover it can be shown that $\widetilde{\mathcal{A}}$ is self-adjoint in $\mathcal{H}^*(\p \widetilde{D}_2)$. 
\begin{lem}
The operator $\widetilde{\mathcal{A}}$ is self-adjoint in $\mathcal{H}^*(\p \widetilde{D}_2)$, {\it i.e.},
$$
( \widetilde\phi,\widetilde{\mathcal{A}}[\widetilde\psi])_{\mathcal{H}^*(\p \widetilde{D}_2)} =
( \widetilde\psi,\widetilde{\mathcal{A}}[\widetilde\phi])_{\mathcal{H}^*(\p \widetilde{D}_2)} 
$$
for $\widetilde\phi,\widetilde\psi \in \mathcal{H}^*(\p \widetilde{D}_2)$.
\end{lem}
\proof
For simplicity, we consider the case when $\widetilde\phi = \varphi_m^c$ and $\widetilde\psi = \varphi_n^c$ only. The other cases can be done similarly.
From \eqref{aaabbbccc}, we have
$
\mathcal{S}_{\widetilde{D}_2}[\varphi_n^c]|_{\p \widetilde{D}_2} = -\frac{\tilde{r}_2}{2n} \varphi_n^c
$. Then, using \eqref{tilde_M_nm} and \eqref{Atilde_matrix}, we have
\begin{align*}
( \varphi_n^c,\widetilde{\mathcal{A}}[\varphi_m^c])_{\mathcal{H}^*(\p \widetilde{D}_2)} &= -( \varphi_n^c,\mathcal{S}_{\p \widetilde{D}_2}\widetilde{\mathcal{A}}[\varphi_m^c])_{-\frac{1}{2},\frac{1}{2}}
\\
&=
-\frac{\tilde{r}_2^{-(n+m-1)}}{8 n m} M_{nm}(\lambda_{D_1} , \widetilde{D}_1).
\end{align*}
So we get the conclusion.
\qed

\smallskip

So $\mathcal{A}$
admits the following spectral decomposition:
$$
\widetilde{\mathcal{A}} = \sum_{n=1}^\infty \lambda_j {\widetilde\psi}_n \otimes {\widetilde\psi}_n
$$
where $\{(\lambda_n, \widetilde\psi_n): n\geq 1\}$ is the set of its eigenvalue-eigenfunction pairs.  We order the eigenvalues in such a way that $|\lambda_j|$ is decreasing and tends to $0$ as $j \to \infty$.  
We remark that all the eigenvalues $\{\lambda_j: j \geq 1\}$ lie in the interval $(-1/2,1/2)$. Moreover, they can be numerically approximated by the eigenvalues of a finite truncation of the infinite matrix $\widetilde{\mathcal{A}}$.

Thanks to \eqref{A_identical_Atilde}, if we let $\psi_n=\widetilde{\psi}_n\circ\Phi$, then we obtain
\be \label{eq-spec-rep}
\mathcal{A} = \sum_{n=1}^\infty \lambda_j {\psi}_n \otimes {\psi}_n.
\ee
It is also worth mentioning that the orthogonality of basis $\{\psi_n\}$ is also preserved.

Using the spectral representation formula (\ref{eq-spec-rep}), we can derive the following result.

\begin{thm} \label{thm1 biosensing}
Assume that Condition \ref{condition0 biosensing} holds and that $D_2$ is in the strong interaction regime, then the scattered field $u^s_{D_2}=u- u_{D_1}$ by the plasmonic particle $D_2$ has the following  representation:
$$
u^s_{D_2} = \mathcal{S}_{D_2,D_1} [\psi],
$$
where $\psi$ satisfies
\beas
\psi = \sum_{j=1}^{\infty}\f{ \left(\nabla u^i(z)\cdot\nu,\psi_j\right)_{\mathcal{H}^*(\p D_2)} \psi_j + O(\delta^2)}{ \lambda_{D_2} - \lambda_j}. 
\eeas
\end{thm}

As a corollary, we obtain the following asymptotic expansion of the scattered field $u-u^i$.
\begin{thm} \label{thm-far field sensing}
The following far field expansion holds:
\beas
(u-u^i)(x) =   \nabla u^i(z)\cdot M(\lambda_{D_1},\lambda_{D_2},D_1,D_2) \nabla G(x,z) +  O\left(\f{\delta^3 }{\textnormal{dist}(\lambda_{D_2},\sigma(\mathcal{A}))}\frac{1}{|x|^2}\right),
\eeas
as $|x|\rightarrow \infty$.
Here, $z$ is the center of mass of $D_2$ and $M(\lambda_{D_1},\lambda_{D_2},D_1,D_2)$ is the polarization tensor satisfying
\be \label{eq-PT_D1D2 sensing}
M(\lambda_{D_1},\lambda_{D_2},D_1,D_2)_{l,m} = \sum_{j=1}^{\infty}\f{(\nu_l,\psi_j)_{\mathcal{H}^*(\p D_2) } (\psi_j,x_m)_{-\f{1}{2},\f{1}{2}} + O(\delta^2)}{ \lambda_{D_2} - \lambda_j },
\ee
for $l,m=1,2$.
\end{thm}
We can introduce the resonant frequency $\omega_j$ for the system of two particles $D_1\cup D_2$ as in Subsection \ref{subsec-plasmonic}.
From the above far field expansion of the scattered field, it is clear that 
when we vary the frequency $\omega$, at certain frequency $\omega$ such that $\lambda_{D_2}(\omega)\approx \lambda_j$ for some $j$ which satisfies the condition that
$$ 
(\nu_l,\psi_j)_{\mathcal{H}^*(\p D_2) } (\psi_j,x_m)_{-\f{1}{2},\f{1}{2}}
\neq 0,
$$
the scattered field will show a sharp peak, which corresponds to the excitation of a plasmonic resonance. Such a frequency is called the (plasmonic) resonant frequency for the system of two particles, which is different from the one for the single plasmonic particle $D_2$. The difference is called the shift of resonant frequency. This shift is due to the interaction of the target particle with the plasmonic particle. As discussed in Subsection \ref{subsec-plasmonic}, the resonant frequencies $\omega_j$ of the two-particle system can also be measured from the far field. They also determines $\lambda_j$ which are eigenvalues of the operator $\mathcal{A}$. 
In the next section, we discuss how to reconstruct the shape of $D_1$ from these recovered eigenvalues.


\section{   { The inverse problem }  } \label{sec-inverse}

In this section, we discuss the inverse problem to  reconstruct the shape of the small unknown particle $D_1$ by using the resonances of the plasmonic particle $D_2$ which interacts with $D_1$. We assume the location of $D_1$ and the permittivity $\epsilon_1$ are known for simplicity. As exlpained in the previous section, we can measure the eigenvalues $\lambda_j$ for $j=1,2,...,J,$ from the far-field measurements. Since the single set of the measurement data is not enough for the reconstruction, we shall make measurements for many different configurations of the two-particles system.  
In Subsection \ref{sec-recon}, we show how the CGPTs of the unknown particle $\widetilde{D}_1$ can be reconstructed from the measurements of $\lambda_j$. In Subsection \ref{sec-optimal}, we explain the optimal control algorithm to recover the shape of $\widetilde{D}_1$
from the CGPTs. 
In this way, we reconstruct the transformed shape $\widetilde{D}_1$ first. Once we find $\widetilde{D}_1$, the original shape of $D_1$ can be easily recovered by using the mapping $\Phi$.
In Subsection \ref{sec-numeric}, we provide several numerical examples.




\subsection{Reconstruction of CGPTs}\label{sec-recon}

In this subsection, we propose an algorithm to reconstruct the CGPTs from measurements of the eigenvalues $\lambda_j$. For ease of presentation, we only consider the first two largest eigenvalues $\lambda_1$ and $\lambda_2$. We denote their measurements by $\mathcal{P}_1$ and $\mathcal{P}_2$, respectively. Note that a single measurement of
$(\mathcal{P}_1,\mathcal{P}_2)$ typically yields very poor reconstruction of the CGPTs due to the lack of information. To overcome this issue, we need to measure the eigenvalues for different configurations of the two particles. Recall the target particle contains the origin. We can rotate it around the origin multiple times and measure $(\mathcal{P}_1,\mathcal{P}_2)$ for each configuration.  
The CGPTs for the target particle after each rotation are related in the following way. 

Define
\begin{align*}
N^{(1)}_{m,n} (\lambda, D) = (M_{m,n}^{cc} - M_{m,n}^{ss}) + i (M_{m,n}^{cs}+M_{m,n}^{sc} ),
\\
N^{(2)}_{m,n} (\lambda, D) = (M_{m,n}^{cc} + M_{m,n}^{ss}) + i (M_{m,n}^{cs}-M_{m,n}^{sc} )
\end{align*}
and let 
$
R_\theta D = \{ e^{i\theta} x : x\in D \},\theta\in[0,2\pi).
$
Then for all integers $m,n$ and all angle parameters $\theta$, we have 
\cite{book3}
\begin{align*}
N^{(1)}_{m,n} ( R_\theta D) = e^{i(n+m)\theta} N^{(1)}_{m,n} (D),
\quad
N^{(2)}_{m,n} (R_\theta D) = e^{i(n-m)\theta} N^{(2)}_{m,n} (D).
\end{align*}

Let us write $d=c\delta$ for some $c>0$. As discussed in subsection \ref{sec-strong-conformal}, if $c$ is of order one, then the deformation of the shape $\widetilde{D}_1$ from $D_1$ is not so strong. So, if the domain $D_1$ is rotated by an angle $\theta$, then the transformed domain will also be rotated by the same amount of angle. So we may (approximately) identify $\widetilde{R_{\theta}D_1} $ with $R_{\theta}\widetilde{D}_1$.

Measuring $\mathcal{P}_j$ for multiple rotation angles $\theta_i$ for $R_\theta\widetilde{D}_1$ will yield a non-linear system of equations that will allow the recovery of the CGPTs associated with $\widetilde{D}_1$.
From the recovered CGPTs, we will reconstruct the shape of $\widetilde{D}_1$. Here, we only consider the shape reconstruction problem. Nevertheless, by using the CGPTs associated with $\widetilde{D}_1$, it is possible to reconstruct the permittivity $\eps_1$ of $\widetilde{D}_1$ in the case it is not a priori given \cite{book3}.

In view of \eqref{Atilde_matrix} and \eqref{Mnm_tilde_decay}, using a standard perturbation method, the asymptotic expansion of the eigenvalue $\lambda_j,j=1,2$, is given by
\begin{align}
\lambda_j =  \lambda_j^0+\lambda_j^1 +\lambda_j^2 + \cdots, \quad \mbox{where } \lambda_j^k = O\big(\left({\tilde{r}_1}/{\tilde{r}_2}\right)^{k+2}\big).
\end{align}
 Each term in the RHS of the above expansion can be computed explicitly. Although we omit the explicit expressions, we mention that they are nonlinear and
depend on CGPTs in the following way: 
\begin{align*}
\lambda_j^0 &= \lambda_j^0(M_{11}),
\\
\lambda_j^1 &= \lambda_j^1(M_{11},M_{12}),
\\
\lambda_j^2 &= \lambda_j^2(M_{11},M_{12},M_{22},M_{13}),
\\
\vdots \ &= \qquad \quad  \vdots
\\
\lambda_j^k &=\ds \lambda_j^k(  \cup_{m+n\leq k+2} \{ M_{mn} \}).
\end{align*}

Suppose we have measurements $\mathcal{P}_1(\theta)$ and $\mathcal{P}_2(\theta)$  for 11 different  rotation angles $\theta_1,\theta_2,...,\theta_{11}$ of the unknown particle $\widetilde{D}_1$. 
We can reconstruct $M_{nm}$ approximately for $m+n\leq 5$. Recall that $M_{mn} = M_{nm}^T$ where subscript $T$ stands for the transpose. 
We look for a set of matrices $\{M_{nm}^{(1)}\}_{m+n\leq 5}$ 
satisfying $[M^{(1)}_{nm}]^T=M^{(1)}_{mn}$ and the the following nonlinear system: for $j=1,2$, 
\beas
\mathcal{P}_{j}(\theta_1) &=&  \sum_{l=0}^{3} \lambda_j^l\Big(      \cup_{m+n\leq l+2} \{ M_{nm}^{(1)}  ({R}_{\theta_1} \widetilde{D}_1) \} \Big),  \\
\mathcal{P}_{j}(\theta_2) &=&  \sum_{l=0}^{3} \lambda_j^l\Big(      \cup_{m+n\leq l+2} \{ M_{nm}^{(1)}  ({R}_{\theta_2} \widetilde{D}_1) \} \Big),
\\
\vdots \quad &=& \qquad\qquad\vdots 
\notag
\\
\mathcal{P}_{j}(\theta_{11}) &=& \sum_{l=0}^{3}  \lambda_j^l\Big(      \cup_{m+n\leq l+2} \{ M_{nm}^{(1)}  ({R}_{\theta_{11}} \widetilde{D}_1) \} \Big).
\eeas
We note that the above equations have 22 independent parameters.
They can be solved by using standard optimization methods.
We expect that 
$$
M_{nm}  = M_{nm}^{(1)} + O\big(\left({\tilde{r}_1}/{\tilde{r}_2}\right)^{6}\big) \quad \mbox{for }m+n\leq 5.
$$

The above scheme can be easily generalized to reconstruct the higher order CGPTs $M_{nm}$. 
This requires more measurement data  $({\mathcal{P}}_1,\mathcal{P}_2)$ from more rotations. Let $k\geq 2$. 
One can see that (using the symmetry $[M^{(k)}_{nm}]^T=M^{(k)}_{mn}$) the set of GPTs $M_{mn}$ satisfying ${m+n\leq 4k+1}$ contains $e_k$ independent parameters, where $e_k$ is given by
\beas
e_k = 16 k^2 + 6k.
\eeas
Therefore, we need $ e_k/2$ pairs of $(\mathcal{P}_1,\mathcal{P}_2)$ to reconstruct these GPTs.
 Let $\{M_{nm}^{(k)}\}_{m+n \leq 4k+1 }$ be the set of matrices satisfying $[M^{(k)}_{nm}]^T=M^{(k)}_{mn}$ and the following system of equations:
\beas
\mathcal{P}_{j}(\theta_i) &=&  \sum_{l=0}^{k-1} \lambda_j^l\Big(      \cup_{m+n\leq l+2} \{ M_{nm}^{(k)}  ({R}_{\theta_i} \widetilde{D}_1) \} \Big), \quad i=1,...,e_k, \ j=1,2.
\eeas
Then we have
$$
M_{nm}  = M_{nm}^{(k)} + O\big(\left({\tilde{r}_1}/{\tilde{r}_2}\right)^{4k+2}\big) \quad \mbox{for }m+n\leq 4k+1.
$$

\subsection{Optimal control approach}\label{sec-optimal}

Now, in order to recover the shape of $\widetilde{D}_1$ from the CGPTs $M_{mn}$, we can minimize the following energy functional 
\be \label{min2c}
 \mathcal{J}_c^{(l)}[B]:= \frac{1}{2} \sum_{H,F\in \{c,s\}}\sum_{n+ m \le k}
 \left| M_{mn}^{HF} (\lambda_{D_1}, B) -  M_{mn}^{HF} (\lambda_{D_1}, \widetilde{D}_1)
 \right|^2 \;,
\ee

We apply the gradient descent method for the minimization. 
We need the shape derivative of the functional $\mathcal{J}_c^{(l)}[B]$.
For $\epsilon$ small, let
$B_\epsilon$ be an $\epsilon$-deformation of $B$, {\it i.e.}, there is a 
scalar function $h \in \mathcal{C}^{1}(\p B)$, such that
 \beas
 \partial B_{\epsilon}:=\{x+\epsilon h(x)\nu(x)\ :  x \in \p B\} . 
 \eeas
According to
\cite{book3,gpt1,gpt2}, the perturbation of the CGPTs due to
the shape deformation is given by
\begin{align}
& M^{HF}_{nm} (\lambda_{D_1}, B_{\epsilon}) -  M^{HF}_{nm} (\lambda_{D_1}, B) 
\nonumber
\\
& \quad  =  \epsilon (k_{\lambda_{D_1}}-1) \int_{ \p B} h(x)\left[\f{\p u}{\p \nu}\Big|_{-} \f{\p v}{\p \nu}\Big|_{-}
 +\frac{1}{k_{\lambda_{D_1}}} \f
 {\p u}{\p T}\Big|_{-}\f{\p v}{\p T}\Big|_{-}\right](x)\, d\sigma(x)+O(\epsilon^2),
\end{align}
where 
\begin{equation} \label{kpbd} k_{\lambda_{D_1}} = (2\lambda_{D_1} +1)/
(2\lambda_{D_1} -1), \end{equation} and $u$ and $v$ are respectively
the solutions to the following transmission  problems:
 \begin{equation}\label{u}
 \left\{
  \begin{array}{ll}
 \Delta u =0 \quad & \mbox{in } \ds B\cup (\R^2 \backslash \overline{B})\;,\\
 \nm
 \ds u|_{+} -u|_{-} =0 \quad &\mbox{on } \p B\;,\\
 \nm
 \ds \f{\p u}{\p \nu}\Big|_{+} -k_{\lambda_{D_1}} \f{\p u}{\p \nu}\Big|_{-} =0 \quad &\mbox{on } \p B\;,\\
 \nm
 \ds (u-H)(x)=O(|x|^{-1}) \quad &\mbox{as } |x|\rightarrow \infty\;,
 \end{array}
\right.
\end{equation}
and
\begin{equation}\label{v}
\left\{
  \begin{array}{ll}
 \Delta v =0 \quad & \mbox{in } \ds B\cup (\R^2 \backslash \overline{B})\;,\\
 \nm
 \ds k_{\lambda_{D_1}} v|_{+} -v|_{-} =0 \quad &\mbox{on } \p B\;,\\
 \nm
 \ds \f{\p v}{\p \nu}\Big|_{+} -\f{\p v}{\p \nu}\Big|_{-} =0 \quad &\mbox{on } \p B\;,\\
 \nm
 \ds (v-F)(x)=O(|x|^{-1}) \quad &\mbox{as } |x|\rightarrow \infty\;.
 \end{array}
\right.
\end{equation}
Here, $\partial /\partial T$ is the tangential derivative. In the case of $M_{nm}^{cs}$,  for example, we put $H=\mbox{Re}\{P_n\} = r^n \cos n\theta$ and $F=\mbox{Im}\{P_m\}=r^n \sin n \theta$. 
The other cases can be handled similarly. 

 Let
\beas
w_{m,n}^{HF}(x) = (k_{\lambda_{D_1}}-1) \left[\f{\p u}{\p \nu}\Big|_{-}
\f{\p v}{\p \nu}\Big|_{-}+\frac{1}{k_{\lambda_{D_1}}}
\f{\p u}{\p T}\Big|_{-}\f{\p v}{\p T}\Big|_{-}\right](x), \quad x \in \p
B\;.
\eeas
The shape derivative of $\mathcal{J}_c^{(l)}$ at $B$ in the direction of $h$
is given by
 \beas
 \langle d_S \mathcal{J}_c^{(l)}[B], h \rangle  =\sum_{H,F\in\{c,s\}} \sum_{m+n \leq k} \delta_{N}^{HF}
 \langle w_{m,n}^{HF}, h \rangle_{L^2(\partial B)} \;,
 \eeas
where
 $$
 \delta^{HF}_{N} = M^{HF}_{nm} (\lambda_{D_1}, B) -  M^{HF}_{nm} (\lambda_{D_1}, \widetilde{D}_1)\;.
 $$

By using the shape derivatives of the CGPTs, we can get an approximation for the matrix $\big(\widetilde{M}_{nm}(\lambda_{D_1}, B_\epsilon)\big)_{n,m=1}^N$ for the slightly deformed shape.  
Next, the shape derivative of $\lambda_j^N(B)$ can be computed by using the standard eigenvalue perturbation theory.
Finally, by applying a gradient descent algorithm, we can minimize, at least locally, the energy functional $\mathcal{J}_c^{(l)}$. Then we get the shape of the original particle $D_1$ using $D_1 = \Phi^{-1}(\widetilde{D}_1)$.

\subsection{Numerical examples} \label{sec-numeric}

In this subsection, we support our theoretical results by numerical examples. 
In the sequel, we set $\delta = 0.001$. We also assume that $B_1$ and $B_2$ are disks of radii $r_1 = \delta $ and $r_2 = 1$, respectively and they are separated by a distance $d=5 \delta$.  Then the ratio $\tilde{r}_1/\tilde{r}_2$ between the transformed radii is approximately $ 0.127$. Note that the ratio is rather small but much larger than the small parameter $\delta$.
We suppose that the material parameter $\varepsilon_1$ of $D_1$ is known and to be given by $\varepsilon_1 = 3$ and so, it holds that $\lambda_{D_1} = 1$. 

We rotate the unknown particle $D_1$ by the angle $\theta_i,i=1,2,...,11$ and get the measurement pair $(\mathcal{P}_1(\theta_i),\mathcal{P}_2(\theta_i))$ for each rotation $\theta_i$, where $\theta_i$ is given by
$$
\theta_i = \frac{2\pi}{11}(i-1), \quad i=1,2,..., 11.
$$ 
We mention that, as discussed in \cite{part1}, we can measure $(\mathcal{P}_1,\mathcal{P}_2)$ from the local peaks of the plasmonic resonant far-field. 
 
\begin{figure*}
\begin{center}
\epsfig{figure=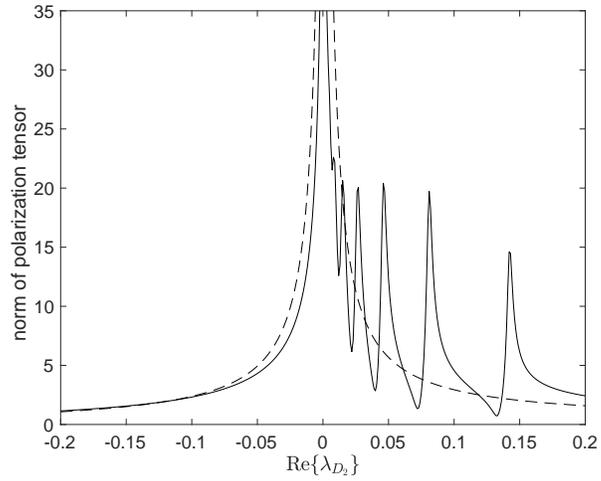,width=9.0cm}
\end{center}
\caption{The magnitude of the polarization tensor.  The dotted line (or solid line) represents the case when the dielectric particle $D_1$ is absent (or presented), respectively.  We set $\mbox{Im}\{\lambda_2\}=0.003$.}
\label{fig:shift}
\end{figure*}

 Figure \ref{fig:shift} shows the shift in the plasmonic resonance.
 In the absence of the dielectric particle $D_1$, the local peak occurs only at  $\lambda_{D_2}=0$. If the particle $D_1$ is presented in a strong regime, then many local peaks appear. By measuring the first two largest values of $\lambda_{D_2}$ at which a local peak appear, we get $(\mathcal{P}_1, \mathcal{P}_2)$ approximately.

\begin{figure*}
\begin{center}
\epsfig{figure=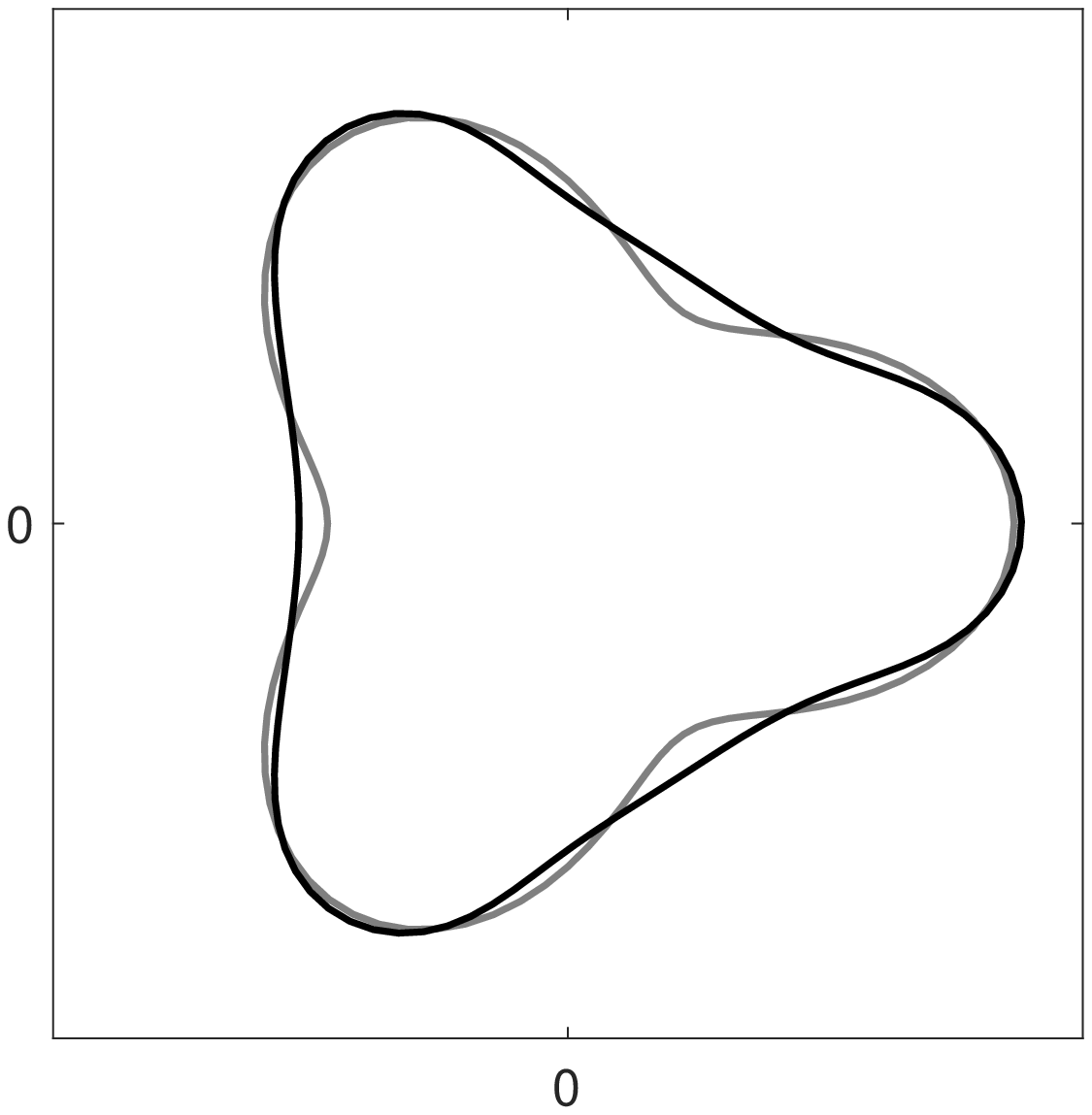,width=5.5cm}
\hskip-0.6cm
\epsfig{figure=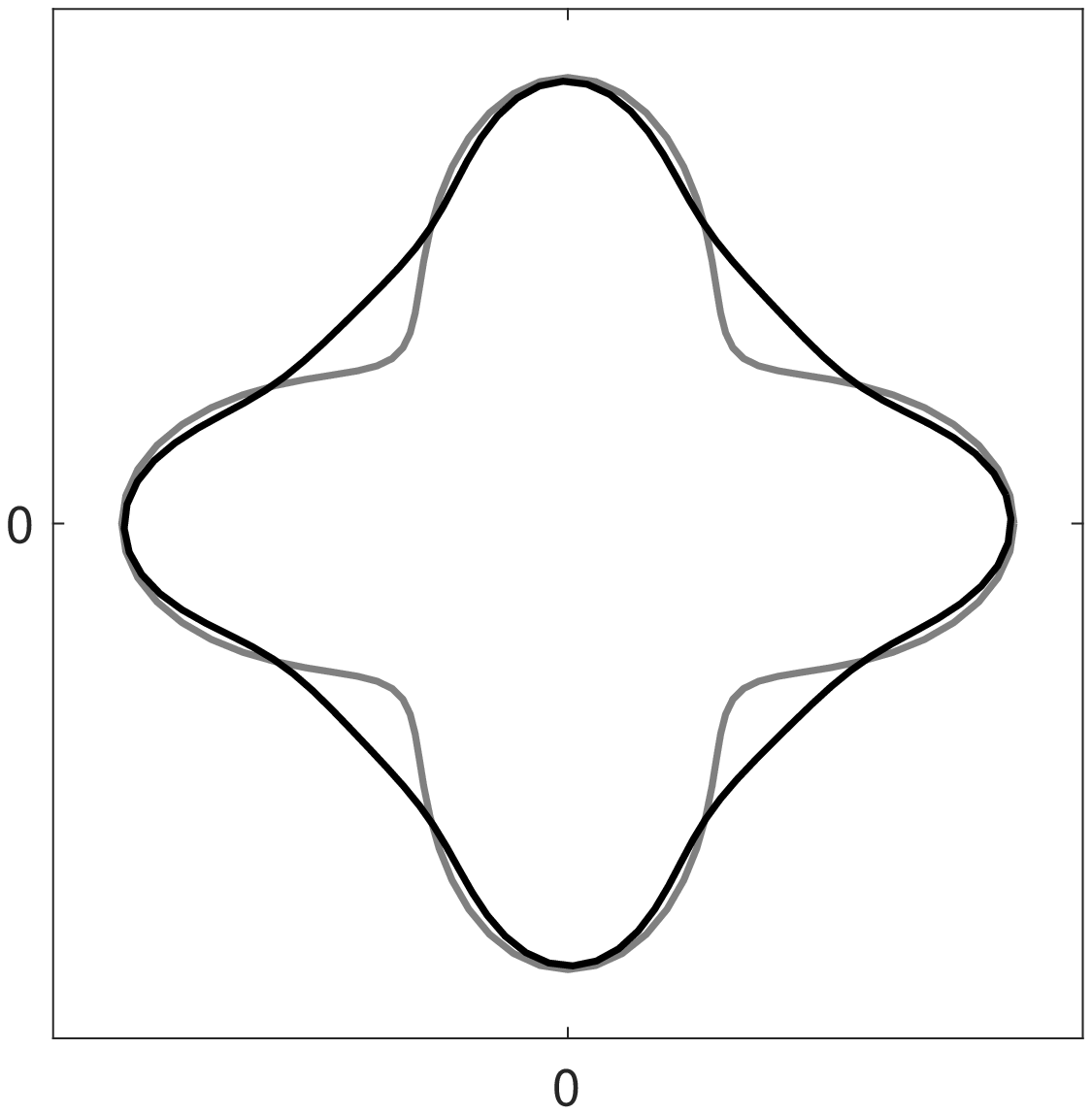,width=5.5cm}
\hskip-0.6cm
\epsfig{figure=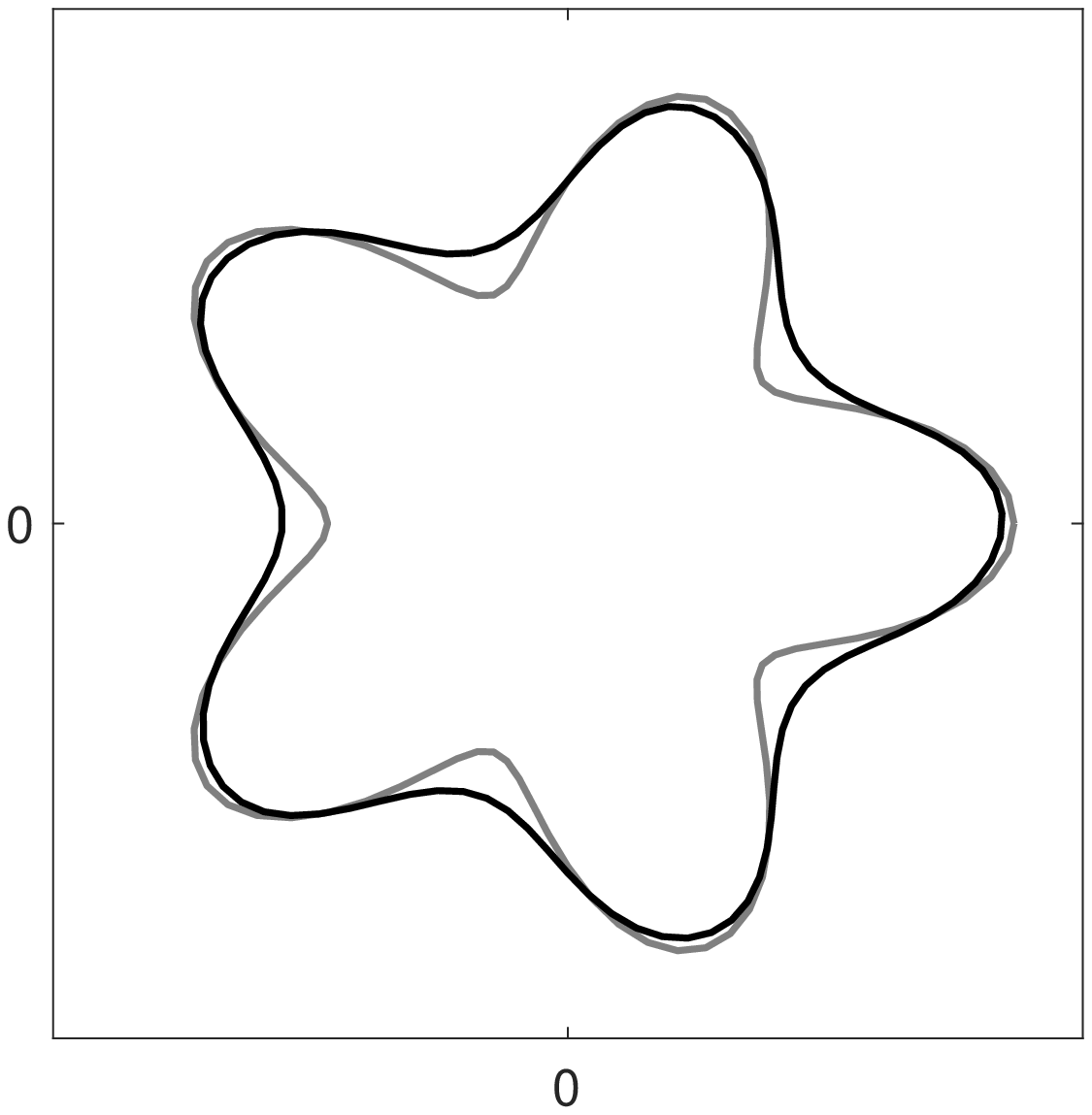,width=5.5cm}
\end{center}
\caption{ Comparison between the original shape (gray) of the particle $D_1$ and the reconstructed one (black). The iteration number is 30.} 
\label{fig:recon}
\end{figure*}

From measurements of $(\mathcal{P}_1, \mathcal{P}_2)$, we recover the contracted GPTs using the algorithm described in subsection \ref{sec-recon}.
We then minimize functional \eqref{min2c} to reconstruct an approximation of $\widetilde{D}_1$. Finally, we use $D_1 = \Phi^{-1}(\widetilde{D}_1)$ to get the shape of $D_1$.
We consider the case of $D_1$ being a flower-shaped particle 
and show comparison between the target shapes and the reconstructed ones, as shown in Figure \ref{fig:recon}. 
We recover the first contracted GPTs up to order 5, {\it i.e.}, $M_{mn}$ for $m+n\leq 5$. 
We take as an initial guess the equivalent ellipse to $\widetilde{D}_1$, determined from the recovered  first order polarization tensor. The required number of iterations is $30$. It is clear that they are in good agreement.


\section{Conclusion}

In this paper, we have made the mathematical foundation of near field sensing complete.  
We have considered the sensing of a small target particle using a plasmonic particle in the strong interaction regime, where the distance between the two particles is comparable to the small size of the target particle. We have introduced a conformal mapping which transforms the two-particle system into a shell-core structure, in which the inner dielectric core corresponds to the target object. Then we have analyzed the shift in the resonance frequencies due to the presence of the inner dielectric core.
We have shown  that this shift  encodes information on the contracted polarization tensors of the core, from which one can reconstruct its shape, and hence the target object. It is worth to mention that although we considered the two dimensional case only in this paper, our conformal mapping approach can be extended to the three dimensional case. Although the Laplacian is not preserved in 3D, there is a nice way to overcome this difficulty \cite{3D-YA}. The extension to the 3D case will be the subject of a forthcoming paper.

\end{document}